\documentclass{elsarticle}

\usepackage{algorithmicx}
\usepackage[ruled]{algorithm}
\usepackage{algcompatible}
\usepackage{algpseudocode}
\usepackage{graphicx}

\usepackage{lineno}
\modulolinenumbers[1]   
\usepackage{amsmath}
\usepackage{amsfonts}
\usepackage{amsbsy}
\usepackage{amscd}
\usepackage{amssymb}
\usepackage{latexsym}
\usepackage[mathscr]{eucal}
\usepackage{indentfirst}

\usepackage{multirow}

\usepackage{url}

\makeatletter \let\cl@part\relax \makeatother

\usepackage[usenames, dvipsnames]{color}
\usepackage{cleveref}

\newtheorem{method}{Method}[section]

\newtheorem{definition}{Definition}[section]

\newtheorem{example}{Example}

\def\CPU{\texttt{CPU}}

\def\method{\texttt{METHOD}}

\def\total{\footnotesize\textnormal{total}}

\def\local{\mbox{\scriptsize local}}

\def\globalmethod{\mbox{\scriptsize glb-meth}}
\def\localmethod{\mbox{\scriptsize  loc-meth}}
\def\localconstruct{\mbox{\scriptsize loc-constr}}
\def\localsolve{\mbox{\scriptsize loc-solv}}

\def\opt{\mbox{\scriptsize opt}}

\journal{Journal of Computational Physics}

\begin{document}

\graphicspath{{figures/}} 

\date{}


\begin{frontmatter}



\title{A local character based method for solving linear systems
       of radiation diffusion problems\tnoteref{label0}}
 \tnotetext[label0]{
                  This work was funded by National Key R\&D Program of China (No. 2017YFA0603903),
                  Science Challenge Project (No. TZ2016002) and
                  National Natural Science Foundation of China (No. 11671051 \& No. 61872380).
                }

 \author[label1]{Shuai Ye}
 \ead{shuaiye09@163.com}
 \author[label2,label3]{Hengbin An\corref{cor1}}
 \ead{an\_hengbin@iapcm.ac.cn}
 \author[label4,label1]{Xinhai Xu}
 \ead{xuxinhai@nudt.edu.cn}
 \cortext[cor1]{Corresponding author. Tel. +86-10-5987 2450.}
 \address[label1]{State Key Laboratory of High Performance Computing, National University of Defense Technology, Changsha 410073, China}
 \address[label2]{Laboratory of Computational Physics,
                  Institute of Applied Physics and Computational Mathematics,
                  Beijing 100094, China}
 \address[label3]{CAEP Software Center for High Performance Numerical Simulation,
                  Beijing 100088, China}
 \address[label4]{Artificial Intelligence Research Center, National Innovation Institute of Defense Technology, Beijing 100071, China.}

\begin{abstract}
  The radiation diffusion problem is a kind of {time-dependent}
  nonlinear equations. For solving the radiation diffusion equations,
  many linear systems are obtained in the nonlinear iterations at each time step.
  The cost of linear equations dominates
  the numerical simulation of radiation diffusion applications,
  such as inertial confinement fusion, etc.
  {Usually,} iterative methods are used to
  solve the linear systems in {a} real application.
  {Moreover, the solution of the previous nonlinear iteration
  or the solution of the previous time step is typically used as
  the initial guess for solving the current linear equations.}
  Because of the strong local character in ICF, with the advancing
  of nonlinear iteration and time step,
  the solution of the linear system changes dramatically in some local domain,
  and changes mildly or even has no change in the rest domain.

  In this paper, a local {character-based} method is proposed
  to solve the linear systems of radiation diffusion problems.
  The proposed method consists of three steps:
  firstly, a local domain (algebraic domain) is constructed;
  secondly, the subsystem on the local domain is solved;
  and lastly, the whole system will be solved.
  Two methods are given to construct the local domain.
  One is based on the spatial gradient{, and} the other is based on the residual.
  Numerical tests for a two-dimensional heat conduction model problem,
  and two real application models, the multi-group radiation diffusion equations
  and the three temperature energy equations, are conducted.
  The test results show {that} the solution time for solving the linear system can be reduced
  dramatically by using the local {character-based} method.
\end{abstract}

\begin{keyword}
  radiation diffusion equations \sep
  linear system \sep
  local character \sep
  spatial gradient  \sep
  residual
\end{keyword}

\end{frontmatter}



\section{Introduction}
\label{sect:introduction}

In the simulation of astrophysics and
inertial confinement fusion (ICF), it is needed to
solve radiation diffusion equations. This is {a
strongly nonlinear system}, and the discretized equations
are solved by some nonlinear iterative methods, such as
Newton kind of methods, Picard method{, and} source iteration method
for multi-group radiation diffusion equations~\cite{baldwin1999iterative, knoll2001nonlinear, mousseau2000physics, yuan2009progress}.
In each nonlinear iteration, some linearized equations should be
solved. In numerical simulations of ICF, usually $\mathcal{O}(10^4)$
time steps are needed to finish the whole simulation, and
a series of linear equations $A x = b$ should be solved.
The solution time for linear equations dominates
the real numerical simulation cost~\cite{ an2009choosing,xu2016adaptive, yuan2009progress}.
It is important to design an efficient method for solving linear equations
in ICF and astrophysics simulations.

Usually, the preconditioned Krylov methods~\cite{saad2003iterative} are used to solve the linear systems
in ICF simulations. For a preconditioned Krylov method, there are four basic procedures:
the construction and application of a preconditioner,
the choice of an initial iterate,
the computation of the next iterate and the check
for stopping criterion~\cite{grinberg2011extrapolation, tromeur2006choice}.
Each procedure has its own effect on the iterative method.
In past years, much research work has been focused on {preconditioning}.
Among them, AMG preconditioner is the most widely used one in the simulation of radiation diffusion equations~\cite{baldwin1999iterative, mousseau2000physics, xu2016adaptive, yuan2009progress}.

For a series of linear systems, several techniques
for choosing an efficient initial iterate
have been proposed in recent years.
For linear systems with the same coefficient matrix and different right-hand sides,
a method for constructing an efficient initial iterate was given in~\cite{chan1997analysis}.
Later, Chan further considered the choice of initial iterate for a series of equations
with different coefficient matrices in~\cite{chan1999galerkin}.
Fischer proposed another method, in which an approximate solution based upon previous solutions was used~\cite{fischer1998projection}. In~\cite{tromeur2006choice},
a model reduction method based on proper orthogonal decomposition (POD) was used to guess a better initial iterate.
The projection methods were also used
in~\cite{barone2009stable,djeddi2017convergence,hall2000proper,markovinovic2006accelerating}.
In~\cite{grinberg2011extrapolation,shterev2015iterative}, a better initial guess was given
by the extrapolation of previous solutions.
An et al. proposed a kind of method for choosing {a} nonlinear initial iterate
for solving {two dimension three-temperature} (2D 3T) heat conduction equations in~\cite{an2009choosing}.
By solving the nonlinear system on a sub-domain where the temperature varied greatly,
the method in~\cite{an2009choosing} can improve the computational efficiency effectively.
In~\cite{gongnonlinear,huang2016nonlinearly,yang2016nonlinear},
some subsystems are constructed to accelerate the convergence
of {nonlinear problems}.

In the simulation of ICF, with time step advancing,
the physical variables,
{such as density,
pressure, temperature, etc.,}
vary greatly in some local domain and
mildly in the rest domain. This is reflected in the radiation diffusion equations,
where the unknown is the temperature or radiation energy.
{From the perspective of solving linear systems
from one nonlinear iteration to the next iteration,
or from one time step to the next, the solution only
varies greatly in some local domain while mildly in the rest.}
In ICF simulations, the initial iterate for solving linear equations is usually set as
the previous nonlinear iterate or previous time step solution.
Therefore, the solution on some local domain where the solution varies greatly has {a}
strong influence on the convergence of an iterative method,
while the solution on the rest domain has little influence on the method.

In this paper, we propose a local character-based method
for solving the linear systems arising from the simulation of radiation diffusion equations.
In this method, a local domain
is first constructed, corresponding to a subset of the global degrees-of-freedom,
on which the solution is expected to vary greatly.
Then, the subsystem on the local domain is solved. At last, the whole system
is solved by using the solution on the local domain.
Two methods are given to construct the local domain.
Specifically, the contributions of this paper are summarized as follows:
\begin{itemize}
\item A local {character-based} method for solving linear systems of radiation diffusion equations is proposed.
\item Two methods are given to construct the local domain.
      One is based on the spatial gradient and the other is based on the residual.
\item Some numerical experiments are conducted
      to verify the effectiveness of the proposed method.
      Compared with the AMG preconditioned Krylov method,
      the local {character-based} method can
      reduce the solution cost in most cases{,} and 40\% cost can be saved at most.
\end{itemize}

The rest of the paper is organized as follows.
The detail of the local {character-based} method is presented in~\cref{sec:lsm},
and also two methods for constructing the local domain are presented in this section.
In~\cref{sec:NR}, some numerical results are given and also the choice of the parameter
for the local {character-based} method is analyzed in this section.
Finally, \cref{sec:summary} is the summary and some remarks of the paper.

\section{The local {character-based} method}
\label{sec:lsm}

In this section, the local {character-based} method is introduced.
The basic algorithm is presented in~\cref{sec:ALS}.
Two methods for constructing the local domain are introduced in the rest.
The first one, which is based on the spatial gradient, will be presented in~\cref{subsect:grad-local}.
The second one, which is based on the residual, will be presented in~\cref{sub:method2}.

\subsection{The local {character-based} algorithm}
\label{sec:ALS}

Consider a linear system resulting from implicit time propagation of discretized
time dependent partial differential equations (PDEs), or the implicit linearization of
a nonlinear system of PDEs:
\begin{eqnarray*}
\label{equ:lin-equs}
A x = b,
\quad
A \in \mathbb{R}^{N \times N},
\quad
x, b \in \mathbb{R}^{N},
\end{eqnarray*}
where
\begin{eqnarray*}
  A =
  \left(
    \begin{array}{cccc}
      a_{11} & a_{12} & \cdots & a_{1N}  \\
      a_{21} & a_{22} & \cdots & a_{2N}  \\
      \vdots   & \vdots & \ddots  & \vdots  \\
      a_{N1} & a_{N2} & \cdots & a_{NN}  \\
    \end{array}
  \right),
  \quad
  x =
  \left(
    \begin{array}{c}
       x_1 \\
       x_2 \\
       \vdots  \\
       x_N \\
    \end{array}
  \right),
  \quad
  b =
  \left(
    \begin{array}{c}
      b_1 \\
      b_2 \\
      \vdots  \\
      b_N \\
    \end{array}
  \right).
\end{eqnarray*}

Eq~(\ref{equ:lin-equs}) is usually obtained from some real applications,
such as the application of radiation diffusion problem, etc.
For solving Eq~(\ref{equ:lin-equs}), an iterative method, particularly a Krylov method is often used.
Because the linear system is obtained in the nonlinear iteration on some time steps,
an effective initial iterate $x^{(0)}$, which is the solution of the last nonlinear iterate
or the solution at last time step, is used.
In the simulation process of ICF,
the temperature and radiation energy only varies greatly in some local domain
and varies mildly in the rest domain.
Therefore, if an iterative method is used to solve Eq~(\ref{equ:lin-equs}),
only some components of the iterate will change greatly
from the initial iterate $x^{(0)}$ to the converged solution.

For purpose of clarity, in the following we introduce some notations.
Let $\Omega = \{ 1, 2, \ldots, N \}$ be the index set of all components (degrees-of-freedom)
of the linear system~(\ref{equ:lin-equs}).
The linear system considered in this paper is discretized from
a partial differential equations, each component corresponds to
one spatial point (for example, a node of the grid) in the computing domain of PDEs.
Therefore, a component will also be called a {\em point} in the following discussion.
$\Omega$ can also be regarded as an algebraic ``domain''.
Assume that, from the initial iterate $x^{(0)}$ to the converged iterate,
the solution varies greatly at $K$ components.
For convenience, the index set of these $K$ components is denoted as
$\Omega_{\local} = \{i_1, i_2, \ldots, i_K \}$.
Each component in $\Omega_{\local}$ is also called a {\em bad point}.
$\Omega_{\local}$ is called a {\em local domain} in the following.
Let $\bar{\Omega}_{\local} = \Omega \setminus \Omega_{\local}$.
Based on the partition
\begin{eqnarray*}
\Omega = \Omega_{\local} \cup \bar{\Omega}_{\local},
\end{eqnarray*}
the linear system~(\ref{equ:lin-equs}) can be partitioned as
\begin{eqnarray}
\label{eq:BEFC}
  \left(
    \begin{array}{cc}
      B & E \\
      F & C
    \end{array}
  \right)
  \left(
    \begin{array}{c}
       x_B \\
       x_C
    \end{array}
  \right)
  =
  \left(
    \begin{array}{c}
       b_B \\
       b_C
    \end{array}
  \right),
\end{eqnarray}
where
\begin{eqnarray*}
x_B =
\left(
\begin{array}{c}
x_{i_1} \\
x_{i_2} \\
\vdots  \\
x_{i_K}
\end{array}
\right),
\quad
b_B =
\left(
\begin{array}{c}
b_{i_1} \\
b_{i_2} \\
\vdots  \\
b_{i_K}
\end{array}
\right),
\end{eqnarray*}
and
\begin{eqnarray*}
\left(
\begin{array}{ll}
B & E \\
F & C
\end{array}
\right)
=
Q^T A Q,
\quad
\left(
\begin{array}{l}
x_B \\
x_C
\end{array}
\right)
=
Q^T x,
\quad
\left(
\begin{array}{l}
b_B \\
b_C
\end{array}
\right)
=
Q^T b.
\end{eqnarray*}
Here $Q \in \mathbb{R}^{n \times n}$ is a permutation matrix
obtained by exchanging the rows $l$ and $i_l$ ($l = 1, 2, \ldots, K$)
of the identity matrix $I \in \mathbb{R}^{n \times n}$.

In Eq~(\ref{eq:BEFC}), $x_B$ is expected to change greatly from
the initial iterate to the converged solution,
while $x_C$ is expected to change mildly from the initial iterate to the converged solution.
In the local character-based method,
first the subsystem
\begin{eqnarray}
\label{equ:loc-equ}
B x_B = b_B - E x_C^{(0)}
\end{eqnarray}
will be solved, and the solution $\bar{x}_B$ will be obtained.
Then
\begin{eqnarray}
\label{equ:sol-assemble}
\tilde{x} =
Q
\left(
\begin{array}{c}
\bar{x}_B \\
x_C^{(0)}
\end{array}
\right)
\end{eqnarray}
will be used as an initial iterate to solve the whole system~(\ref{equ:lin-equs}) iteratively.
The specific local character-based algorithm is described by Algorithm~\ref{alg:local-character}.

\begin{algorithm}[htb]
\caption{Local {character-based} algorithm}
\label{alg:local-character}
\begin{algorithmic}[1]
\State{\textbf{Input}: $A$, $x^{(0)}$ and $b$}.
\State{Construct the subset $\Omega_{\local}$.}
\State{Construct the permutation matrix $Q$.}
\State{Construct the sub-matrix $B$, $C$, $E$ and $F$.
       Partition the right hand side into $b_B$ and $b_C$,
       and partition the initial iterate into $x_B^{(0)}$ and $x_C^{(0)}$.}
\State{Solve the subsystem~(\ref{equ:loc-equ}), obtain the solution $\bar{x}_B$.}
\State{Assemble $\tilde{x}$ defined in (\ref{equ:sol-assemble}) by using $\bar{x}_B$ and $x_C^{(0)}$.}
\State{{Apply several Gauss-Seidel iterations to
       Equations (\ref{equ:lin-equs}) so that $\tilde{x}$ is updated.}}
\If {{$\tilde{x}$ satisfies the convergence criterion}}
    \State {{Let $x=\tilde{x}$.}}
\Else
    \State{{Solve the whole system~(\ref{equ:lin-equs}) by using the initial iterate $\tilde{x}$,
       and obtain the solution $x$.}}
\EndIf
\State{\textbf{Output}: $x$.}
\end{algorithmic}
\end{algorithm}

In Algorithm~\ref{alg:local-character}, there are two key ingredients:
the first one is the construction of the subset $\Omega_{\local}$ in Line 2;
the second one is the solution of the local subsystem in Line 5.
The construction of $\Omega_{\local}$ is the base of the algorithm.
{By $\Omega_{\local}$,} the whole system can be partitioned into two subsystems,
and particularly the local subsystem is defined. By solving the local subsystem,
a solution on $\Omega_{\local}$ is obtained.
After assembling the whole solution $\tilde{x}$ as the initial iterate,
the whole system will be solved by an iterative method.
Because the solution on $\bar{\Omega}_{\local}$ is expected to change mildly,
$x_C^{(0)}$ is expected to be very near to the final solution.
Therefore $\tilde{x}$ will be very near to the solution.
Actually, in some of our numerical cases,
the residual corresponding to $\tilde{x}$ is so small that
the convergence criteria is satisfied,
and $\tilde{x}$ is a final solution.

The construction of the local domain, or the subset $\Omega_{\local}$,
plays a key role in Algorithm~\ref{alg:local-character}.
In the following, two methods will be introduced to construct the local domain.

\subsection{{Gradient-based} local domain construction}
\label{subsect:grad-local}

The linear system~(\ref{equ:lin-equs}) is discretized
from a partial differential equations (such as radiation diffusion equations).
For the solution of a partial differential equations,
its gradient can be defined. For example, if $T^n(x,y)$ is a solution
of a two dimensional radiation diffusion equation at time step $n$,
then the gradient of $T^n(x,y)$ is defined as
$\nabla T^n(x,y) = \left( \partial T^n / \partial x, \partial T^n / \partial y \right)$.
The norm of the gradient, $\|\nabla T^n(x,y)\|$,
can be used to predict the variation of the solution from time step $n$ to time step $n+1$.

\begin{figure}[h]
 \begin{center}
 \includegraphics[height=250pt]{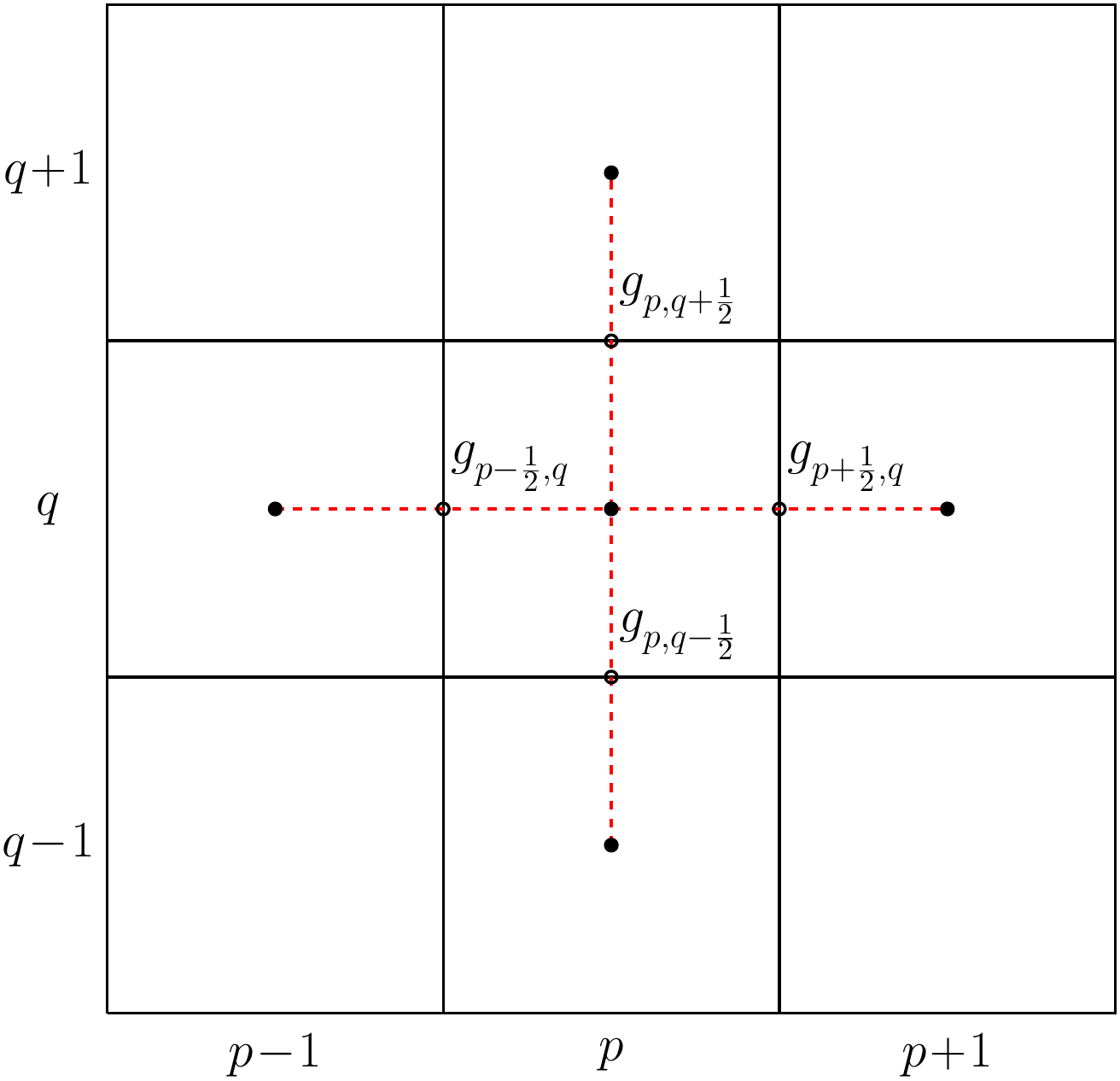}
 \caption{\label{fig:five-scheme}
          Gradient for five-point scheme.
         }
 \end{center}
\end{figure}

For the linear system~(\ref{equ:lin-equs}),
the similar norm of the gradient can be defined,
which can be used to construct the local domain $\Omega_{\local}$.
For a point $i$, if the norm of the gradient of the initial iterate
is large at this point, then it is likely that the solution may change greatly around this point,
which will be considered as an element for the local domain $\Omega_{\local}$.
Specifically, for a given row $i$ of the matrix $A$,
the non-zero entries of this row reflect the adjacency relationship
of the point $i$ to other points.
Therefore, the gradient $g_i$ at point $i$ (component $i$) is defined by
\begin{eqnarray}
\label{equ:grad}
g_i =
\sum_{\substack{1 \le j \le N \\  a_{ij} \neq 0}}
{\left| x^{(0)}_i - x^{(0)}_j \right|}.
\end{eqnarray}
Note that this definition is not the {exact} gradient of a field
(such as temperature, density, etc.)
in {two-dimensional or three-dimensional} space, but just a mimicry.
{More precisely}, this is very similar to
{the} $l_1$ norm of the gradient at each point $i$.
For example, if $A$ is obtained by the {five-point} difference scheme,
then the summation in~(\ref{equ:grad}) includes four terms, as shown in Fig.~\ref{fig:five-scheme},
in which $p$ and $q$ are used as indices for cells.
{By definition} in (\ref{equ:grad}), the gradient {of} cell $(p,q)$ is
\begin{eqnarray*}
g_{p,q}
&=& |g_{p-\frac{1}{2},q}| + |g_{p+\frac{1}{2},q}| + |g_{p,q-\frac{1}{2}}| + |g_{p,q+\frac{1}{2}}| \\
&=& |x_{p-1,q} - x_{p,q}| + |x_{p+1,q} - x_{p,q}| + |x_{p,q-1} - x_{p,q}| + |x_{p,q+1} - x_{p,q}|.
\end{eqnarray*}

Based on the definition of the gradient~(\ref{equ:grad}), the local domain
can be constructed. For point $i$, if
\begin{eqnarray*}
g_i > \alpha \max_{1 \leq j \leq N}{g_j},
\end{eqnarray*}
then $i$ will be a member of the local domain index set $\Omega_{\local}$.
Here $\alpha$ is a prescribed parameter.
The specific algorithm is shown in Algorithm~\ref{alg:gradient-local-domain}.

\begin{algorithm}[htbp]
\caption{{Gradient-based} local domain construction}
\label{alg:gradient-local-domain}
\begin{algorithmic}[1]
\State{\textbf{Input}: $A$, $x^{(0)}$, $\alpha$.}
\State{\textbf{Initialize}: $g_{\max} = 0$, $\Omega_{\local} = \{ \, \}$, $K=0$.}
\For {$i = 1:N$}
   \State{Compute $g_i$ by equation~(\ref{equ:grad}).}
   \If {$g_i > g_{\max}$}
       \State{$g_{\max} = g_i$.}
   \EndIf
\EndFor
\For {$i = 1:N$}
   \If {$g_i > \alpha g_{\max}$}
       \State{$\Omega_{\local} = \{\Omega_{\local}, i\}$.}
       \State{$K=K+1$.}
   \EndIf
\EndFor
\State{\textbf{Output}: $\Omega_{\local}$, $K$.}
\end{algorithmic}
\end{algorithm}

\subsection{{Residual-based} local domain construction}
\label{sub:method2}

The second method for constructing the local domain is based on
the residual of the linear system. The method consists of two steps.
First find the components with large residual,
and the bad points set, or the local domain, $\Omega_{\local}$, will be constructed.
Second consider the influence of the bad points on other
degrees-of-freedom, and the set of bad points then extended to
include other degrees-of-freedom
on which the solution is also likely to change.

\subsubsection{Finding bad points}
\label{sec:badpoint}

To iteratively solve the linear system~(\ref{equ:lin-equs}),
a commonly used convergence criteria is given by
\begin{eqnarray}
\label{equ:rel-res-conv-criteria}
\|b - A x^{(k)}\|_2
\le
\epsilon \|b\|_2,
\end{eqnarray}
where $\epsilon$ is a prescribed tolerance and
 $x^{(k)}$ is an iterative solution.
Since
\begin{eqnarray*}
\|b - A x^{(k)}\|_2 \le \sqrt{N} \|b - A x^{(k)}\|_{\infty},
\end{eqnarray*}
the condition~(\ref{equ:rel-res-conv-criteria}) will be satisfied if
\begin{eqnarray*}
\|b - A x^{(k)}\|_{\infty}
\le
\frac{\epsilon \|b\|_2}{\sqrt{N}},
\end{eqnarray*}
or
\begin{eqnarray*}
\max_i \left| b_i - \sum_{j=1}^N a_{i,j} x_j^{(k)} \right|
\le
\frac{\epsilon \|b\|_2}{\sqrt{N}},
\end{eqnarray*}
which is equivalent to
\begin{eqnarray*}
\left| b_i - \sum_{j=1}^N a_{i,j} x_j^{(k)} \right|
\le
\frac{\epsilon \|b\|_2}{\sqrt{N}},
\quad
i = 1, 2, \ldots, N.
\end{eqnarray*}

Now, by the initial iterate $x^{(0)}$ and the convergence criteria tolerance $\epsilon$,
a component $i$ is defined as a {\em bad point} if
\begin{eqnarray*}
\left| b_i - \sum_{j=1}^N a_{i,j} x_j^{(0)} \right|
>
\frac{\epsilon \|b\|_2}{\sqrt{N}}.
\end{eqnarray*}
Specifically, the initial local domain $\Omega_{\local}$ is defined by
\begin{eqnarray}
\label{equ:resi-local-domain-0}
\Omega_{\local} =
\left\{
i : i = 1, 2, \ldots, N;
\left| b_i - \sum_{j=1}^N a_{i,j} x_j^{(0)} \right|
>
\frac{\epsilon \|b\|_2}{\sqrt{N}}
\right\}.
\end{eqnarray}

\subsubsection{Expanding the local domain}
\label{subsubsec:expand-domain}

When the local subsystem of the bad points is solved,
the solution at the bad points is changed.
This then has further influence on points adjacent to the bad points.
Therefore, the set of the bad points will further be extended.
That is, the subset of the bad points (local domain) obtained
by~(\ref{equ:resi-local-domain-0}) will further be expanded.

As discussed in Subsection~\ref{sec:ALS}, assume that $\bar{x}_B$ is
the solution of the local domain subsystem, and let $\tilde{x}$ defined by~(\ref{equ:sol-assemble}).
The $j$-th residual corresponding to $x^{(0)}$ and $\tilde{x}$ are given, respectively, as
\begin{eqnarray*}
r^{(0)}_j = b_j - \sum_{i=1}^N a_{j,i} x_i^{(0)},
\quad
\mbox{and}
\quad
\tilde{r}_j = b_j - \sum_{i=1}^N a_{j,i} \tilde{x}_i.
\end{eqnarray*}
By the residuals, the influence of the local domain to the point $j$ can be given.
\begin{definition}
\label{def:influence}
Assume that $\Omega_{\local}$ is a local domain.
For one point $j \in \bar{\Omega}_{\local}$,
let $r_j^{(0)}$ be the residual at point $j$ with the initial iterate,
and let $\tilde{r}_j$ be the residual at point $j$ after solving the local domain subsystem.
If $\left| {r}_j^{(0)} \right| \le \frac{\epsilon \|b\|_2}{\sqrt{N}}$ and
$\left| \tilde{r}_j \right| > \frac{\epsilon \|b\|_2}{\sqrt{N}}$,
then the local domain $\Omega_{\local}$ has influence on the point $j$;
otherwise, the local domain $\Omega_{\local}$ has no influence on the point $j$.
\end{definition}

It is easy to see that
\begin{eqnarray*}
\tilde{r}_j
&=& b_j - \sum_{i=1}^N a_{j,i} \tilde{x}_i \\
&=& b_j - \sum_{l=1}^K a_{j,i_l} \tilde{x}_{i_l} - \sum_{l=K+1}^N a_{j,i_l} \tilde{x}_{i_l} \\
&=& b_j - \sum_{l=1}^K a_{j,i_l} \bar{x}_{l} - \sum_{l=K+1}^N a_{j,i_l} x_{i_l}^{(0)} \\
&=& b_j - \sum_{l=1}^K a_{j,i_l} \left( \bar{x}_{l} - x_{i_l}^{(0)} \right) - \sum_{l=1}^N a_{j,i_l} x_{i_l}^{(0)} \\
&=& r_j^0 - \sum_{l=1}^K a_{j,i_l} \left( \bar{x}_{l} - x_{i_l}^{(0)} \right).
\end{eqnarray*}
Assume that
\begin{eqnarray*}
\left| \bar{x}_{l} - x_{i_l}^{(0)} \right| \le \left| x_{i_l}^{(0)} \right|,
\quad
l = 1, 2, \ldots, K,
\end{eqnarray*}
which means the relative variation of the solution at all the bad points is less than one
(This is satisfied in our numerical tests).
Then it is easy to see that
\begin{eqnarray}
\label{equ:extend-local-condition}
\left| \tilde{r}_j \right| \le \left| {r}_j^{(0)} \right| +  \sum_{l=1}^K \left| a_{j,i_l} x_{i_l}^{(0)} \right|,
\quad
j = 1, 2, \ldots, N.
\end{eqnarray}
This inequality shows that if
\begin{eqnarray*}
\left| {r}_j^{(0)} \right| + \sum_{l=1}^K \left| a_{j,i_l} x_{i_l}^{(0)} \right|
\le
\frac{\epsilon \|b\|_2}{\sqrt{N}},
\end{eqnarray*}
then the local domain will have not influence on point $j$
(by Definition~\ref{def:influence}).
Otherwise, the local domain will be expected to have influence on point $j$.
In this case, the local domain will be extended.
Specifically, assume that $j \in \bar{\Omega}_{\local}$ is a neighbor of $\Omega_{\local}$.
If
\begin{eqnarray*}
\left| {r}_j^{(0)} \right| +  \sum_{l=1}^K \left| a_{j,i_l} x_{i_l}^{(0)} \right|
>
\frac{\epsilon \|b\|_2}{\sqrt{N}},
\end{eqnarray*}
then $j$ will be incorporated into the local domain subset $\Omega_{\local}$.

Note that the extension process can be implemented several times
until the size of the subset $\Omega_{\local}$ no longer increases.
In practice, a few times of expanding is enough. The specific process
is given in Algorithm~\ref{alg:residual-local-domain}.

In Algorithm~\ref{alg:residual-local-domain},
two steps are used to construct the local domain $\Omega_{\local}$.
The first step is to construct $\Omega_{\local}$ by the component residuals,
and those components with large residuals are incorporated into $\Omega_{\local}$.
The second step is to extend $\Omega_{\local}$ by at most $E_{\max}$ iterations.
In each expanding iteration, some neighbors of $\Omega_{\local}$ are further
incorporated into $\Omega_{\local}$.

\begin{algorithm}[htbp]
\caption{Residual based local domain construction}
\label{alg:residual-local-domain}
\begin{algorithmic}[1]
\State{\textbf{Input}: $N$, $A$, $b$, $x^{(0)}$, $\epsilon$, $E_{\max}$.}
\State{\textbf{Initialize}: $\Omega_{\local} = \{ \, \}$, $\mathcal{N} = \{ \, \}$, $K=0$.}
\For {$i = 1:N$}
   \State{Compute the initial residual $r_i^{(0)} = b_i - \sum_{j=1}^N a_{i,j} x_j^{(0)}$.}
   \If {$|r_i^{(0)}| > \frac{\epsilon \|b\|_2}{\sqrt{N}}$}
       \State{$\Omega_{\local} = \{\Omega_{\local}, i\}$.}
       \State{$K = K+1$.}
   \EndIf
\EndFor
\For {$m = 1:E_{\max}$}
  \State{Let $\mathcal{N}$ be the neighbors of $\Omega_{\local}$.}
  \State{Let $\mathcal{N}_{\local} = \{ \, \}$, $K_0 = 0$.}
  \ForAll {$j \in \mathcal{N}$}
     \If {$\left| {r}_j^{(0)} \right| +  \sum_{l=1}^K \left| a_{j,i_l} x_{i_l}^{(0)} \right| > \frac{\epsilon \|b\|_2}{\sqrt{N}}$}
       \State{$\mathcal{N}_{\local} = \{\mathcal{N}_{\local}, j\}$.}
       \State{$K_0 = K_0+1$.}
     \EndIf
  \EndFor
  \If {$K_0 > 0$}
     \State{$\Omega_{\local} = \{ \Omega_{\local}, \mathcal{N}_{\local} \}$.}
     \State{$K = K + K_0$.}
  \Else
     \State{break;}
  \EndIf
\EndFor
\State{\textbf{Output}: $\Omega_{\local}$, $K$.}
\end{algorithmic}
\end{algorithm}

\begin{example}
\label{exam:res-local-character}
The following example shows the specific implementation of Algorithm~\ref{alg:residual-local-domain}.
Consider a $9 \times 9$ linear system with the coefficient matrix given by
\begin{eqnarray*}
A =
\left(
\begin{matrix}
1 & -\frac{1}{2} &  &  &   &  &  &  &              \\
-\frac{1}{2} & 1 & -\frac{1}{3} & & & & & &        \\
  & -\frac{1}{3} & 1 & -\frac{1}{4} & & & & &      \\
  &  & -\frac{1}{4} & 1 & -\frac{1}{5} & & & &     \\
  &  &  & -\frac{1}{5} & 1  & -\frac{1}{6} & & &   \\
  &  &  &  & -\frac{1}{6} & 1 & -\frac{1}{7} & &   \\
  &  &  &  &  & -\frac{1}{7} & 1 & -\frac{1}{8} &  \\
  &  &  &  &  &  & -\frac{1}{8} & 1 & -\frac{1}{9} \\
  & &  &  &   &  &  & -\frac{1}{9}  & 1
\end{matrix}
\right)
\end{eqnarray*}
The solution and the initial iterate are set, respectively, as
\begin{eqnarray*}
x =
\left(
\begin{matrix}
1e{-1} \\
1e{-2} \\
1e{-3} \\
1e{-4} \\
1e{-5} \\
1e{-6} \\
1e{-7} \\
1e{-8} \\
1e{-9}
\end{matrix}
\right),
\quad
x^{(0)} =
\left(
\begin{matrix}
1           \\
1e{-1}      \\
1.001e{-3}  \\
1.001e{-4}  \\
1.001e{-5}  \\
1.001e{-6}  \\
1.001e{-7}  \\
1.001e{-8}  \\
1.001e{-9}
\end{matrix}
\right).
\end{eqnarray*}
The right hand side is set to $b = Ax$.

Consider the application of Algorithm~\ref{alg:residual-local-domain}
to this linear system.
The convergence criterion is set to $\epsilon=1e{-5}$,
and the parameter for maximal number of extension is set to $E_{\max} = 6$.
It is easy to check that ${\epsilon \|b\|_2}/\sqrt{N} = 3.44e{-7}$.
The initial residual is
\begin{eqnarray*}
r^{(0)} =
b - A x^{(0)} =
\left(
\begin{matrix}
-8.55e{-1}   \\
 3.66e{-1}   \\
 3.00e{-2}   \\
 1.52e{-7}   \\
 1.02e{-8}   \\
 6.81e{-10}  \\
 4.41e{-11}  \\
 2.61e{-12}  \\
 1.11e{-13}
\end{matrix}
\right)
\end{eqnarray*}

In the first step of Algorithm~\ref{alg:residual-local-domain},
points 1, 2 and 3 are marked as the bad points, i.e.,
$\Omega_{local}=\{1,2,3\}$.
In the second step of Algorithm~\ref{alg:residual-local-domain},
$\Omega_{local}$ is further extended.
The loop of extension process is implemented 4 times.
After each loop, the main results are listed in the following table.
In the fourth iteration, $\mathcal{N} = \{7\}$ and
$|r^{(0)}_7|+|a_{76} \cdot x^{(0)}_6| = 1.430e{-7} < {\epsilon \|b\|_2}/\sqrt{N} = 3.44e{-7}$,
therefore, $K_0 = 0$ and no point will be incorporated into $\Omega_{\local}$.
The loop will break.

\begin{table}[h]
\setlength{\tabcolsep}{4.0mm}
\caption{
\label{table:omega-local-extend}
The extension process of $\Omega_{\local}$.}
\begin{center}
\begin{tabular}{cccccl}
\hline
$m$   & $\mathcal{N}$ & $\mathcal{N}_{\local}$ & $\Omega_{\local}$ & $K$ &
$| {r}_j^{(0)} | +  \sum_{l=1}^K | a_{j,i_l} x_{i_l}^{(0)} |$ \\
\hline
1 & $\{ 4 \}$ & $\{ 4 \}$ & $\{ 1, 2, 3, 4 \}$       & 4 & $|r^{(0)}_4|+|a_{43} \cdot x^{(0)}_3| = 2.504e{-4}$ \\
\hline
2 & $\{ 5 \}$ & $\{ 5 \}$ & $\{ 1, 2, 3, 4, 5 \}$    & 5 & $|r^{(0)}_5|+|a_{54} \cdot x^{(0)}_4| = 2.003e{-5}$ \\
\hline
3 & $\{ 6 \}$ & $\{ 6 \}$ & $\{ 1, 2, 3, 4, 5, 6 \}$ & 6 & $|r^{(0)}_6|+|a_{65} \cdot x^{(0)}_5| = 1.669e{-6}$ \\
\hline
4 & $\{ 7 \}$ & $\emptyset$ & $\{ 1, 2, 3, 4, 5, 6 \}$  & 6 & $|r^{(0)}_7|+|a_{76} \cdot x^{(0)}_6| = 1.430e{-7}$ \\
\hline
\end{tabular}
\end{center}
\end{table}

\end{example}

\subsection{{Parallelization of the local character-based methods}}
\label{subsect:para}

For parallel computing, assume that $p$ processors are used.
For simplicity, assume that $N = \bar{N} \times p$, where $\bar{N}$
is an integer.
The matrix $A$ is partitioned by rows, that is
\begin{eqnarray*}
A =
\left[
A_{\underline{1}}^{T} \quad
A_{\underline{2}}^{T} \quad
\cdots \quad
A_{\underline{p}}^{T}
\right]^{T},
\end{eqnarray*}
where $A_{\underline{l}}$ is a $\bar{N} \times N$ matrix, $l = 1, 2, \ldots, p$.
The solution and right-hand side are partitioned correspondingly
\begin{eqnarray*}
x =
\left[
x_{\underline{1}}^{T} \quad
x_{\underline{2}}^{T} \quad
\cdots \quad
x_{\underline{p}}^T
\right]^{T},
\quad
b =
\left[
b_{\underline{1}}^{T} \quad
b_{\underline{2}}^{T} \quad
\cdots \quad
b_{\underline{p}}^T
\right]^{T},
\end{eqnarray*}
where $x_{\underline{l}}$ and $b_{\underline{l}}$
are $\bar{N}$ vectors, $l = 1, 2, \ldots, p$.
For convenience, let
\begin{eqnarray*}
\mathcal{I}_{\underline{l}} = \left\{ (l-1)\bar{N}, \ldots, l \bar{N} \right\}
\end{eqnarray*}
be the component index set of the $l$-th processor, $l = 1, 2, \ldots, p$.

The key work in the parallel version of the algorithm is to construct the local domain
$\Omega_{\local}$. When the local domain is constructed, the local domain subsystem
can be obtained easily. For solving the local domain subsystem,
any parallel solution method can be used.
In the following, Algorithm~\ref{alg:gradient-local-domain-para} and
Algorithm~\ref{alg:residual-local-domain-para} are given for the gradient based
and residual based local domain construction in parallel computing case.

\begin{algorithm}[htbp]
\caption{Gradient-based local domain construction in parallel}
\label{alg:gradient-local-domain-para}
\begin{algorithmic}[1]
\State{\textbf{Input}: $A_{\underline{l}}$, $x_{\underline{l}}^{(0)}$, $\mathcal{I}_{\underline{l}}$, $\alpha$.}
\State{\textbf{Initialize}: $g_{\max} = 0$, $\Omega_{\local} = \{ \, \}$, $K=0$.}
\For {$i \in \mathcal{I}_{\underline{l}}$}
   \State{Communicate all adjacency components of $x_i^{(0)}$ from other processors.} \label{line:grad-meth-comm-neigh}
   \State{Compute $g_i$ by equation~(\ref{equ:grad}).}
   \If {$g_i > g_{\max}$}
       \State{$g_{\max} = g_i$.}
   \EndIf
\EndFor
\State{Reduction $g_{\max}$ by max.} \label{line:grad-meth-comm-reduc1}
\For {$i \in \mathcal{I}_{\underline{l}}$}
   \If {$g_i > \alpha g_{\max}$}
       \State{$\Omega_{\local} = \{\Omega_{\local}, i\}$.}
       \State{$K=K+1$.}
   \EndIf
\EndFor
\State{Reduction $\Omega_{\local}$ and $K$ by sum.} \label{line:grad-meth-comm-reduc2}
\State{\textbf{Output}: $\Omega_{\local}$, $K$.}
\end{algorithmic}
\end{algorithm}
Compared to Algorithm~\ref{alg:gradient-local-domain},
some communications should be implemented for parallel case,
which is shown in Line~\ref{line:grad-meth-comm-neigh}, \ref{line:grad-meth-comm-reduc1}, and \ref{line:grad-meth-comm-reduc2} in Algorithm~\ref{alg:gradient-local-domain-para}.

\begin{algorithm}[htbp]
\caption{Residual based local domain construction in parallel}
\label{alg:residual-local-domain-para}
\begin{algorithmic}[1]
\State{\textbf{Input}: $N$, $\bar{N}$, $A_{\underline{l}}$, $b_{\underline{l}}$, $x^{(0)}_{\underline{l}}$, $\mathcal{I}_{\underline{l}}$, $\epsilon$, $E_{\max}$.}
\State{\textbf{Initialize}: $\Omega_{\local} = \{ \, \}$, $\mathcal{N} = \{ \, \}$, $K=0$.}
\For {$i \in \mathcal{I}_{\underline{l}}$}
   \State{Communicate all adjacency components of $x_i^{(0)}$ from other processors.} \label{line:resi-meth-comm-neigh1}
   \State{Compute the initial residual $r_i^{(0)} = b_i - \sum_{j=1}^N a_{i,j} x_j^{(0)}$.}
   \If {$|r_i^{(0)}| > \frac{\epsilon \|b\|_2}{\sqrt{N}}$}
       \State{$\Omega_{\local} = \{\Omega_{\local}, i\}$.}
       \State{$K = K+1$.}
   \EndIf
\EndFor
\State{Reduction $\Omega_{\local}$ and $K$ by sum.}
\label{line:resi-meth-comm-reduc1}
\For {$m = 1:E_{\max}$}
  \State{Let $\mathcal{N} = \{ \mbox{the neighbors of } \Omega_{\local} \} \cap \mathcal{I}_{\underline{l}}$.}
  \State{Let $\mathcal{N}_{\local} = \{ \, \}$, $K_0 = 0$.}
  \ForAll {$j \in \mathcal{N}$}
     \State{Communicate all adjacency components of $x_j^{(0)}$ from other processors.}\label{line:resi-meth-comm-neigh2}
     \If {$\left| {r}_j^{(0)} \right| +  \sum_{l=1}^K \left| a_{j,i_l} x_{i_l}^{(0)} \right| > \frac{\epsilon \|b\|_2}{\sqrt{N}}$}
       \State{$\mathcal{N}_{\local} = \{\mathcal{N}_{\local}, j\}$.}
       \State{$K_0 = K_0+1$.}
     \EndIf
  \EndFor
  \If {$K_0 > 0$}
     \State{$\Omega_{\local} = \{ \Omega_{\local}, \mathcal{N}_{\local} \}$.}
     \State{$K = K + K_0$.}
     \State{Reduction $\Omega_{\local}$ and $K$ by sum.}
     \label{line:resi-meth-comm-reduc2}
  \Else
     \State{break;}
  \EndIf
\EndFor
\State{\textbf{Output}: $\Omega_{\local}$, $K$.}
\end{algorithmic}
\end{algorithm}

Compared to Algorithm~\ref{alg:residual-local-domain},
it is needed to have some communication
in Line~\ref{line:resi-meth-comm-neigh1}, \ref{line:resi-meth-comm-reduc1}, \ref{line:resi-meth-comm-neigh2}, and \ref{line:resi-meth-comm-reduc2} in Algorithm~\ref{alg:residual-local-domain-para}.

It should be noted that for both Algorithm~\ref{alg:gradient-local-domain-para} and Algorithm~\ref{alg:residual-local-domain-para}, some global reductions are needed to be implemented.
For example, the computation of the right hand side norm, etc.

\section{Numerical Results}
\label{sec:NR}

In this section, one test model and two suites of real application linear systems
are used to test the effectiveness of the local {character-based} method.
The test model is a two-dimensional nonlinear heat conduction equation.
The two suites of real application linear systems
are respectively the multi-group radiation diffusion equations
and the three temperature energy equations.

All experiments are carried out on an E5-2620 CPU,
which is clocked at 2.10 GHz and has a 16 GB memory.
E5-2620 has a total of 12 cores on two sockets.

The local {character-based} method is compared with the typical solution method.
In the test, the compared method is a BoomerAMG preconditioned GMRES method.
For BoomerAMG, Falgout coarsening (a combination of CLJP and the classical Ruge-Stuben coarsening)
and ``classical'' interpolation~\cite{ruge1987algebraic} are used,
the maximal number of levels is set to 8,
and one symmetric Gauss-Seidel iteration is used for both pre-relaxation and post-relaxation. The maximal GMRES iteration number is set to 80,
and the Krylov dimension is 40.

In the local {character-based} method, the local domain subsystem is solved by
the BoomerAMG preconditioned GMRES method with completely the same parameters
for solving the global system. One Gauss-Seidel iteration is applied to smooth $\tilde{x}$.
In the local {character-based} methods,
the parameters $\alpha$ in Algorithm~\ref{alg:gradient-local-domain} and
$E_{\max}$ in Algorithm~\ref{alg:residual-local-domain} are predefined by users,
and some numerical analysis is given for different $\alpha$ and $E_{\max}$.

For the local {character-based} method,
after {solving} the sub-domain system,
the assembled solution $\tilde{x}$ in (\ref{equ:sol-assemble}) satisfies
the convergence criteria in most of the tests.
Therefore, no further iteration is implemented for the global system in most cases.

{For the purpose} of presenting the results,
the following notations are used to represent the specific methods:
\begin{itemize}
\item $\method0$: BoomerAMG preconditioned GMRES method.
\item $\method1$: Local {character-based} method in which the local domain is constructed by Algorithm~\ref{alg:gradient-local-domain}.
\item $\method2$: Local {character-based} method in which the local domain is constructed by Algorithm~\ref{alg:residual-local-domain}.
\end{itemize}
Furthermore, the following notations will be used to report the results.
\begin{itemize}
\item $N$: {the} scale of the global linear system.
\item $K$: {the} scale of the sub-domain linear system in {the} local {character-based} method.
\item $\eta = \frac{K}{N} \times 100\%$: percentage of the scale of the sub-domain linear system to the scale of the global system.
\item $\CPU_{\total}^{\globalmethod}$: CPU time for solving the linear system by BoomerAMG preconditioned GMRES method.
\item $\CPU_{\total}^{\localmethod}$: CPU time for solving the linear system by the local {character-based} method.
\item $\CPU_{\localconstruct}$: CPU time for constructing the local domain in {the} local {character-based} method.
\item $\CPU_{\localsolve}$: CPU time for solving the sub-domain system in {the} local {character-based} method.
\item $S_{\CPU}$: {the} speedup of CPU time for solving the linear system by
      the local {character-based} method to BoomerAMG preconditioned GMRES method.
\end{itemize}

\subsection{Test of {two-dimensional} heat conduction equation}
\label{subsec:2DheatT}

Consider the following {two-dimensional} (2-D) heat conduction equations:
\begin{eqnarray*}
\label{eq:2DheatT}
\left\{
\begin{array}{ll}
\dfrac{\partial T}{\partial t} =
\nabla \cdot \left( \kappa(T) \nabla T \right), &  0 < x < 1, \; 0 < y < 1, \; 0 < t < \infty \\
T(0,x,y) = T_0(x,y),  &  0 < x < 1, \;  0 < y <1, \; t = 0 \\
T(t,0,y) = T_l, & x = 0, \; 0 < t < \infty \\
T(t,1,y) = T_r, & x = 1, \; 0 < t < \infty \\
\dfrac{\partial T}{\partial y} = 0, & y = 0, \; y = 1, \; 0 < t < \infty
\end{array}
\right.
\end{eqnarray*}
{In this test},
$T_l$, $T_r$, $\kappa(T)$ and $T_0(x, y)$ are set respectively as

$$
T_l = 1,
\quad
T_r = 1 \times 10^{-4},
\quad
\kappa(T) = T^{3.5},
\quad
T_0(x,y) = e^{-100x} \times T_l + T_r.
$$

By using the backward Euler method to discretize in time, one obtains

\begin{eqnarray}
\label{eqn:2D-heat-temp-dist}
\dfrac{T^{n+1} - T^n}{\Delta t} =
\left(
\nabla \cdot \left( \kappa(T) \nabla T \right)
\right)^{n+1}
\end{eqnarray}
The mesh is $99 \times 99$,
and a five-point finite difference method is used to discretize {in space},
and a nonlinear system is obtained.

In this test, the time step size is $\Delta t = 10^{-2}$.
The simulation is implemented from physical time 0 to 1,
and the total number of time steps is 100.
The Picard method is used to solve the system,
and the linear system in Picard iteration
will be solved by $\method0$, $\method1$, or $\method2$.
The parameter $\alpha$ in $\method1$ is set to $10^{-4}$,
{and} the $E_{\max}$ in $\method2$ is set to 1.
For clarity, the time advancing and nonlinear iteration process is given
in Algorithm~\ref{alg:2d-cond-sol}.

The $n$ loop in Algorithm~\ref{alg:2d-cond-sol} is for time advancing,
and the $s$ loop is for nonlinear iteration.
Eq.~(\ref{eqn:picard-equs}) is the Picard linearized system.
The convergence tolerance for the linear system is $\epsilon = 10^{-10}$, and
the stopping criterion for the nonlinear Picard iteration is $\left\|T^{n+1,s+1}-T^{n+1,s}\right\|_2<10^{-8}$.
In the computation, the previous nonlinear iterate is used as
the initial guess for solving the current linear system.

\begin{algorithm}[htbp]
\caption{Solution for 2-D heat conduction equation}
\label{alg:2d-cond-sol}
\begin{algorithmic}[1]
\State{\textbf{Input}: The number of mesh for $x$ and $y$ directions, $N_x$, $N_y$,
       time step size $\Delta t$ and the initial $T^0$.}
\For {$n = 0, 1, \ldots$}
   \State{$T^{n+1,0} = T^n$.}
   \For {$s = 0, 1, \ldots$}
     \State{Obtain $T^{n+1,s+1}$ by solving the following linear equations}
     \begin{eqnarray}
     \label{eqn:picard-equs}
       \dfrac{T^{n+1,s+1} - T^n}{\Delta t} =
       \nabla \cdot \left( \kappa(T^{n+1,s}) \nabla T^{n+1,s+1} \right)
     \end{eqnarray}
     \If {(nonlinear iteration converge)}
       \State{break;}
     \EndIf
     \State{$T^{n+1,s} = T^{n+1,s+1}$.}
   \EndFor
\EndFor
\end{algorithmic}
\end{algorithm}

\begin{figure}[htbp]
  \centering
  \includegraphics[width=1\linewidth]{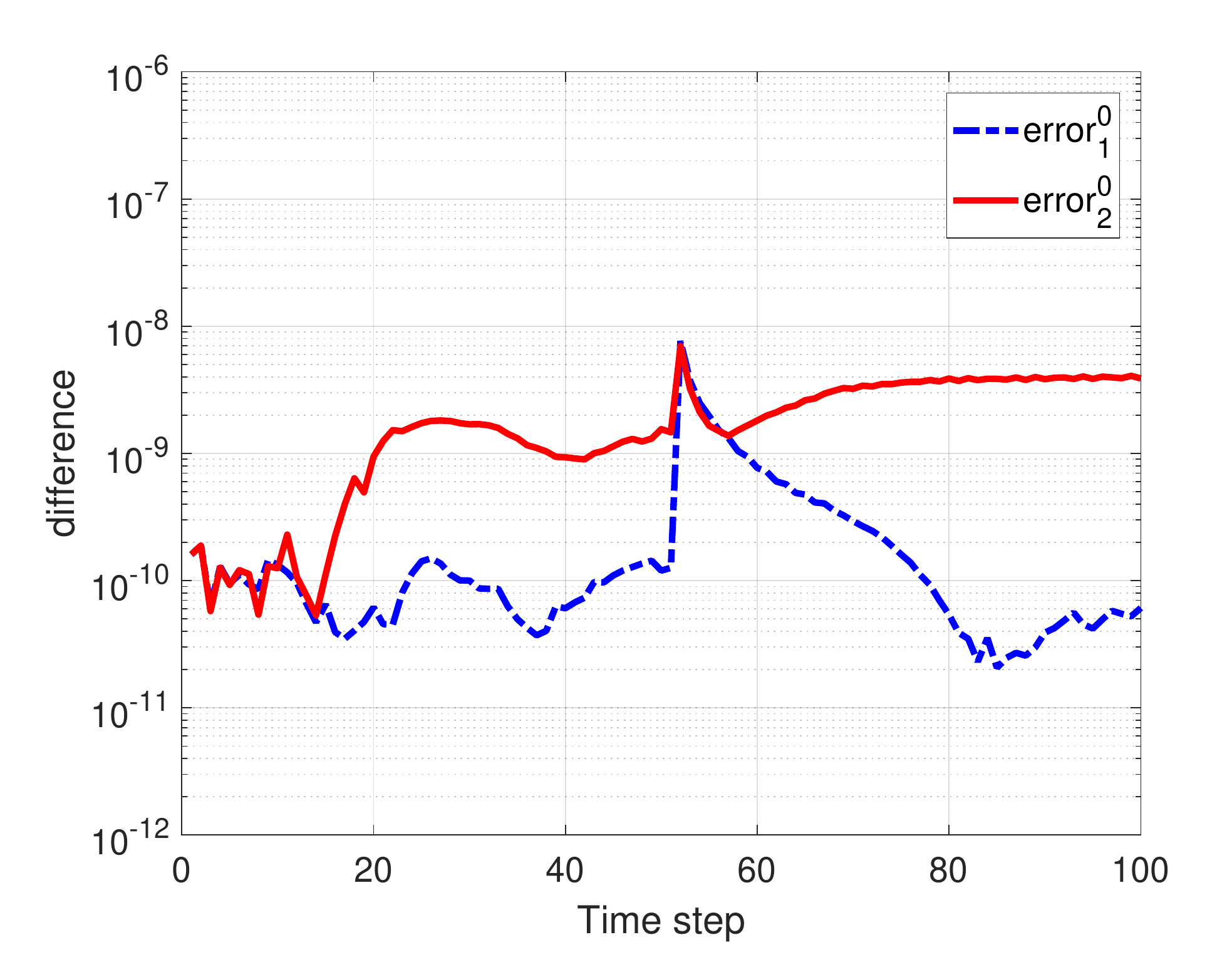}
  \caption{
  The differences of the solution obtained by $\method0$, $\method1$, and $\method2$ for two dimensional heat conduction  problem.}
  \label{fig:error-2D}
\end{figure}

First {a} comparison between the solutions obtained by the three different methods is conducted.
Fig.~\ref{fig:error-2D} shows
the difference in solution obtained by different methods as a function of time, given by
\begin{eqnarray*}
e_1^0 = \left\| T_{\method1} - T_{\method0}\right\|_2 / \left\| T_{\method0} \right\|_{2},
\quad
e_2^0 = \left\| T_{\method2} - T_{\method0}\right\|_2 / \left\| T_{\method0} \right\|_{2}.
\end{eqnarray*}
Here, $T_{\method0}$, $T_{\method1}$, and $T_{\method2}$
are the solutions obtained by $\method0$, $\method1$,
and $\method2$, respectively, at each time step.
From this figure, one can see that for all 100 time steps,
the maximal value of $e_1^0$ and $e_2^0$ are, respectively, $7.6\times 10^{-9}$ and $7.0\times10^{-9}$.
This shows that the approximation solution obtained by different methods are very close
and they can be considered almost the same.

In the test, if $\method0$, $\method1$, or $\method2$ is used
to solve all the linear systems for all 100 time steps,
the total CPU time for solving the linear systems
is respectively $39.19s$, $28.09s$, and $29.76s$.
Therefore, the speedups of
$\method1$ and $\method2$ to $\method0$ are about 1.40 and 1.32, respectively.

\begin{figure}[htbp]
  \centering
  \includegraphics[width=0.9\linewidth]{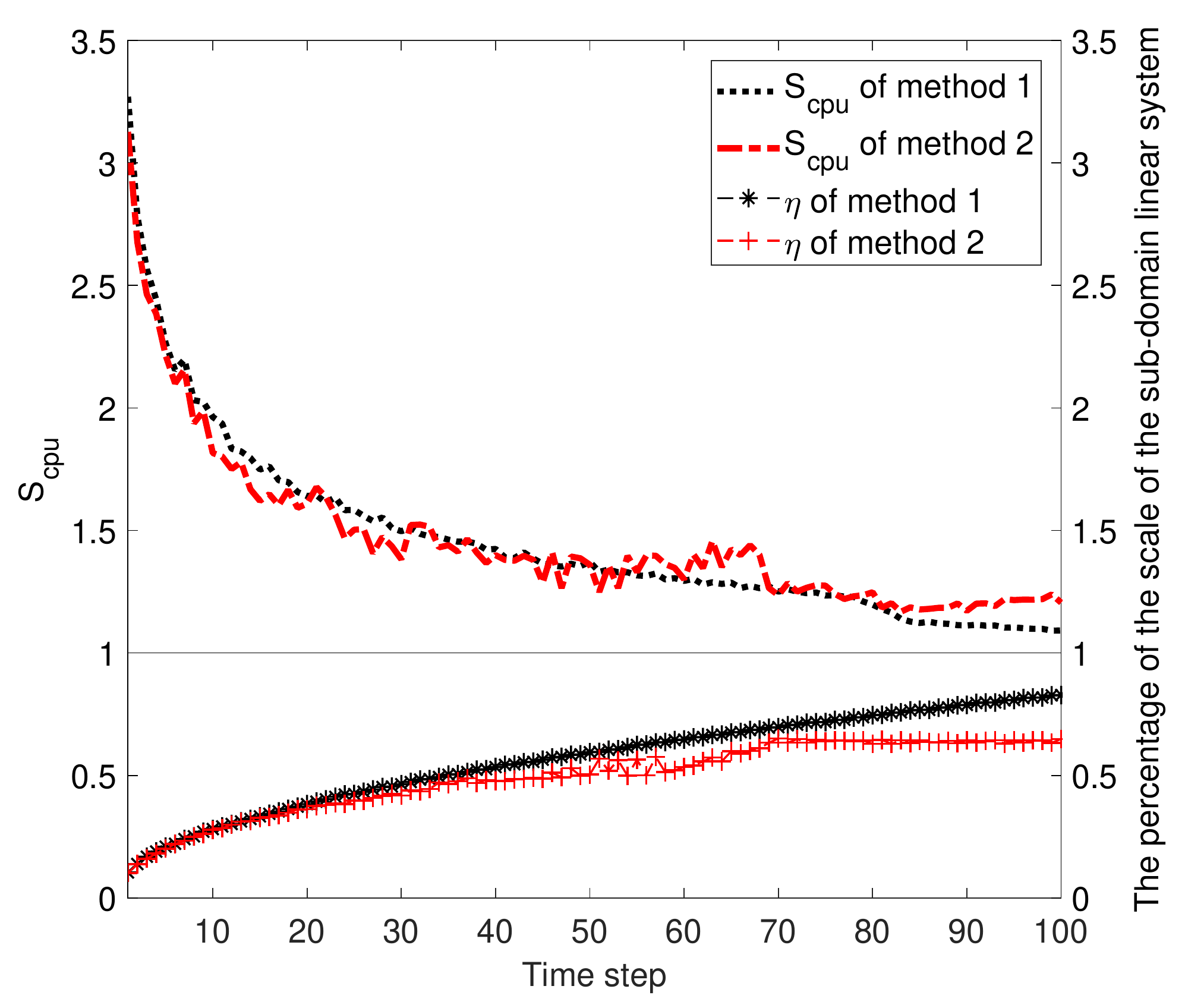}
  \caption{The average speedup and average percentage of the scale of the sub-domain linear system to the scale of the global system for two local character-based methods.}
  \label{fig:2DnonlinearT}
\end{figure}

Fig.~\ref{fig:2DnonlinearT} shows the average speedups of the two local character-based methods
at each time step. At same time, the average percentage of the scale of the sub-domain
linear system to the scale of the global system, $\eta$,
for two local character-based methods is shown in the figure.
From this figure, one can see that both methods have the similar performance in the simulation.
With time advancing, the speedup deceases, and $\eta$ increases.
At later simulation time, $\method2$ performs a little better than $\method1$.

\begin{figure}[htbp]
  \centering
  \begin{minipage}[t]{0.9\linewidth}
  \includegraphics[width=0.9\linewidth]{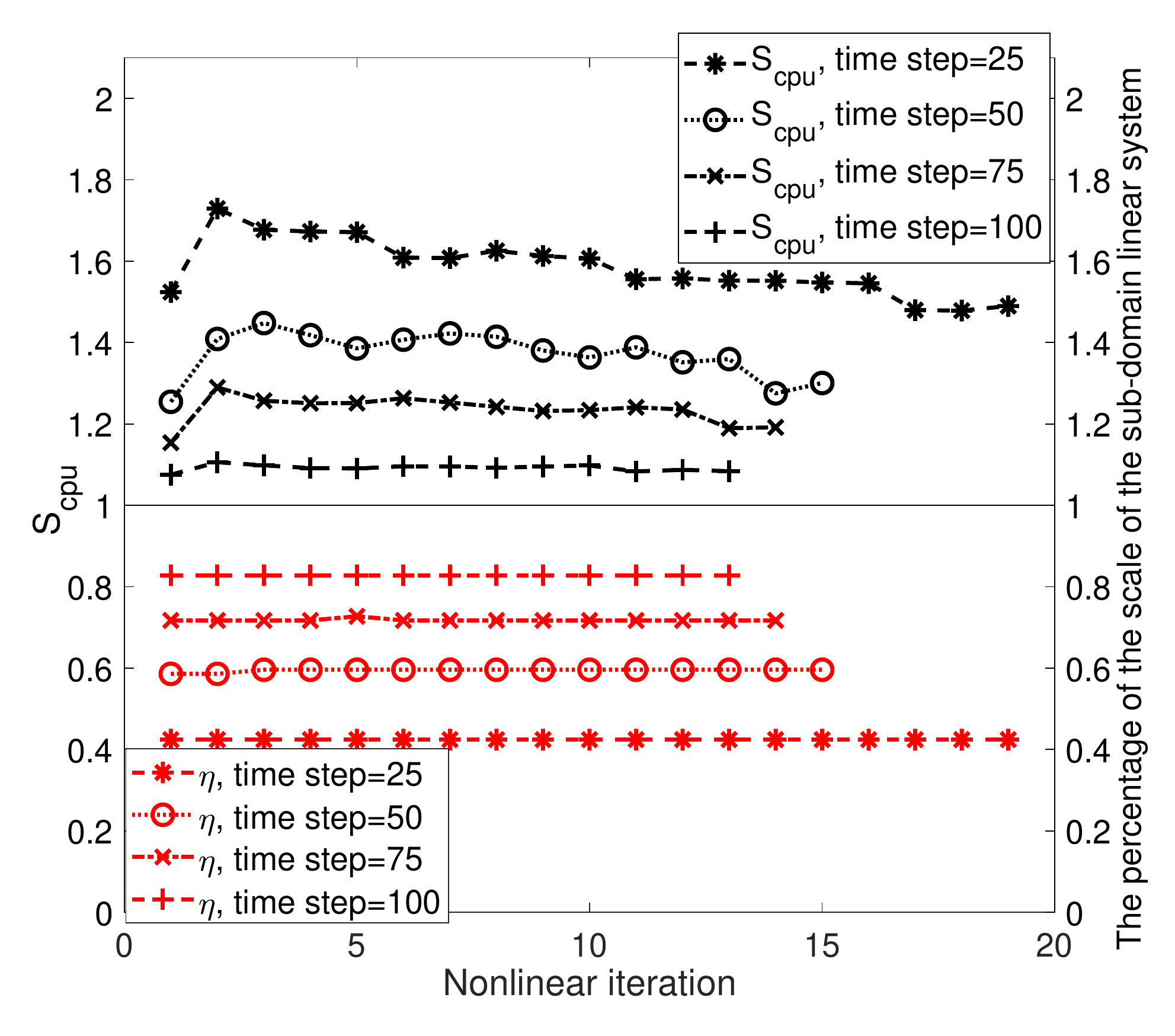}
  \end{minipage}
  \begin{minipage}[t]{0.9\linewidth}
  \includegraphics[width=0.9\linewidth]{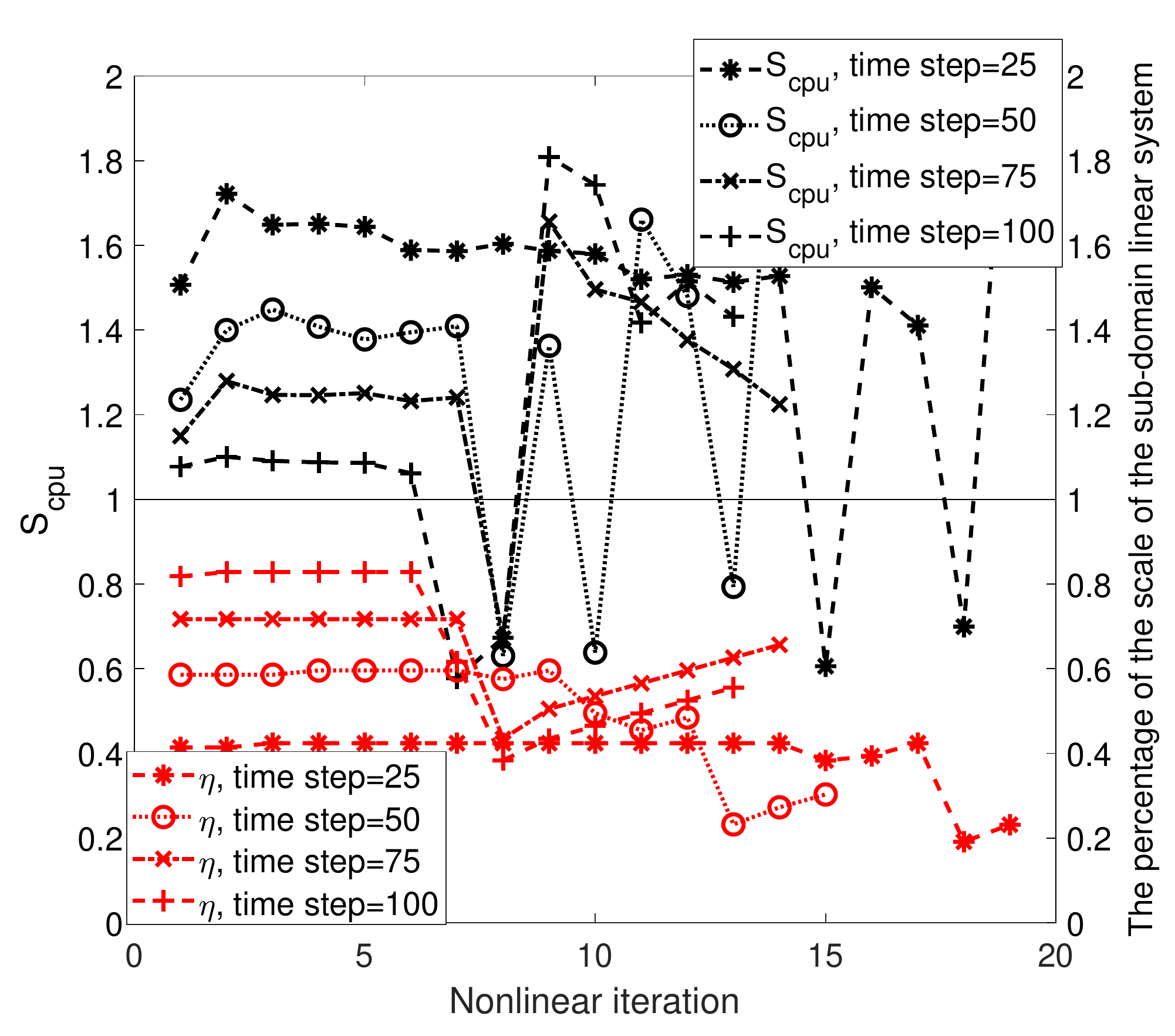}
  \end{minipage}
  \caption{The percentage of the scale of the subsystem and
  the speedup of the local {character-based} methods
  for solving linear systems of the 2-D heat conduction equations at time steps 25, 50, 75, and 100.
  The upper is the result of $\method1$, the lower is the result of $\method2$.}
  \label{fig:2Dlinear}
\end{figure}

Fig.~\ref{fig:2Dlinear} shows the specific speedup and
percentage of the scale of the sub-domain linear system to
the scale of the global system at each nonlinear iteration
for the two local character-based methods.
In the figure, the results at four time steps are given.
From this figure, one can see that $\method1$ have a relatively stable performance
at each time step. At the same time, the percentage of the sub-domain linear system
is stable for all nonlinear iterations at each time.
With time advancing, the performance of $\method1$ decreases.
For $\method2$, one can see that it has similar performance
for the early 6 nonlinear iterations at all time steps.
For the later nonlinear iterations, $\method2$ performs unstable.
At about three nonlinear iterations, its performance deteriorates.

\begin{figure}[htbp]
  \centering
  \includegraphics[width=1.1\linewidth]{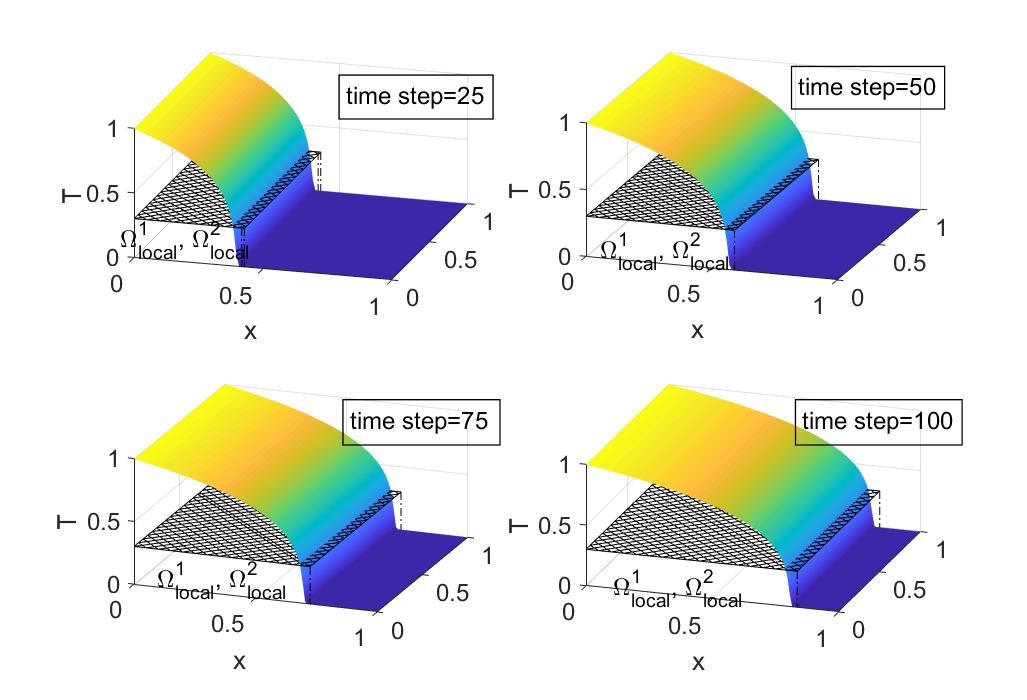}
  \caption{The local domain for two dimensional heat conduction problem at time steps 25, 50, 75, and 100 (the first nonlinear iteration).}
  \label{fig:local-domain-2D}
\end{figure}

For 2-D heat conduction equation, {the local domains} determined
by $\method1$ and $\method2$ are almost the same.
As examples, Fig.~\ref{fig:local-domain-2D} shows the specific local domain in the first nonlinear iteration
at time steps 25, 50, 75, and 100.
At the same time, the solutions at these two time steps are shown in the figure.
One see can that the local domain $\Omega_{\local}^1$ and $\Omega_{\local}^2$ are almost the same,
and they just contain the region where the solution varies greatly.
Here $\Omega_{\local}^1$ represents the local domain obtained by Algorithm~\ref{alg:gradient-local-domain}
and $\Omega_{\local}^2$ represents the local domain obtained by Algorithm~\ref{alg:residual-local-domain}.

\subsection{{Scalability  of the local character-based methods for the 2-D heat conduction equation}}
\label{subsect:scalability}

In this section, the scalability test of the local character-based methods
are done for the 2-D heat conduction equation.
In the test, the scale on each processor is fixed as $512 \times 512$,
and the number of processors varies from 1 to 256.
For all tests, the time step size is fixed as $\Delta t = 10^{-2}$,
and the number of time steps is 100. The test results are listed in Table~\ref{table:par-2D}.
In the table, the following notations are used:
\begin{itemize}
\item $\CPU_{\total}$: CPU time for solving all the linear systems in the simulation.
\item $N_{\mbox{\footnotesize{eq}}}$: Total number of linear equations in the simulation.
\item $\overline{\CPU_{\total}}$ : The average CPU time for solving one linear system.
\item $S_{\CPU}$: The speedup of $\method1$ and $\method2$ to $\method0$ (based on the average CPU time for solving one linear system).
\end{itemize}

\begin{table}
\setlength{\tabcolsep}{4.0mm}
\caption{
\label{table:par-2D}
Scalability test of $\method0$, $\method1$, and $\method2$ for 2-D heat conduction equation}
\begin{center}
\begin{tabular}{ccccccc}
\hline
 \#cores & Method     &  $N$    & $\CPU_{\total}$ & $N_{\mbox{\footnotesize{eq}}}$ & $\overline{\CPU_{\total}}$ & $S_{\CPU}$ \\
\hline
        & $\method0$ &          & 1087.45 & 2757 & 0.39 & -    \\
   1    & $\method1$ & $512^2$  &  693.38 & 2757 & 0.25 & 1.56 \\
        & $\method2$ &          &  760.71 & 2758 & 0.28 & 1.42 \\
\hline
        & $\method0$ &          & 2364.74 & 3281 & 0.72 & -    \\
    4   & $\method1$ & $1024^2$ &  927.42 & 3288 & 0.28 & 2.55 \\
        & $\method2$ &          & 1148.46 & 3319 & 0.35 & 2.06 \\
\hline
        & $\method0$ &          & 3783.12 & 3655 & 1.04 & -    \\
    16  & $\method1$ & $2048^2$ & 1105.13 & 3770 & 0.29 & 3.42 \\
        & $\method2$ &          & 1608.81 & 3844 & 0.42 & 2.35 \\
\hline
        & $\method0$ &          & 4172.75 & 3955 & 1.06 & -    \\
   64   & $\method1$ & $4096^2$ & 1100.34 & 4253 & 0.26 & 3.79 \\
        & $\method2$ &          & 1789.89 & 4293 & 0.42 & 2.33 \\
\hline
        & $\method0$ &          & 4621.16 & 4052 & 1.14 & -    \\
   256  & $\method1$ & $8192^2$ & 1283.76 & 4741 & 0.27 & 3.60 \\
        & $\method2$ &          & 2173.07 & 4676 & 0.46 & 2.13 \\
\hline
 \end{tabular}
\end{center}
\end{table}

From Table~\ref{table:par-2D}, one can see that $\method1$ and $\method2$ scale better than $\method0$
by comparing the the total CPU time for solving all the linear equations, $\CPU_{\total}$.
By comparing the numbers of linear equations for different methods, we find that
more linear equations are solved for $\method1$ and $\method2$ than $\method0$.
By comparing the average CPU time for solving one linear equation,
one can see that $\method1$ and $\method2$ perform much better than $\method0$.
Correspondingly, the speedups of $\method1$ and $\method2$ (for solving one linear equation)
are from 1.42 to 3.79. With the increase of the scale, the speedup increases for both
$\method1$ and $\method2$.

\subsection{Test of multi-group radiation diffusion equations}
\label{sub:result1}

In this subsection, the local {character-based} method will
be used to solve the multi-group radiation diffusion equations
in real applications. {The linear equations for the test} are extracted from
the radiation hydrodynamics code---Lared-integration,
which is used for inertial confinement fusion (ICF) study~\cite{song2015lared}.

\subsubsection{Multi-group radiation diffusion equations}

The multi-group radiation diffusion equations are {a set of $G$ coupled} partial differential equations, which are given as follows:
\begin{eqnarray}
\label{eq:MG}
\frac{\partial E_g}{\partial t} = \nabla \cdot \kappa_g \nabla E_g + 4 \pi \sigma_g B_g - c \sigma_g E_g,
\quad g = 1,\cdots, G
\end{eqnarray}
In the equations,
\begin{itemize}
\item $G$ is the number of {energy} groups.
\item $E_g$ is the radiation energy density of the $g$-th group.
\item $\kappa_g$ is the radiation diffusion coefficient of the $g$-th group.
\item $B_g$ is the integration of the Planck function of the $g$-th group.
\item $\sigma_g$ is the cross section of the $g$-th group.
\end{itemize}

In real applications, the multi-group radiation diffusion equations are coupled to electron-ion equations,
and they form a strong nonlinear system.
For more details about the radiation diffusion model,
see~\cite{hang2013convergence, yuan2009progress}.
We will focus on the solution for {the} multi-group radiation diffusion equations in this paper.
Equations~(\ref{eq:MG}) are a time-dependent nonlinear system.
The backward Euler method is used to discretize in time.
For solving the nonlinear system after temporal discretization,
a nonlinear iteration (based on Picard-Newton linearization) is used,
and in each nonlinear iteration, the source iteration method is used
to solve the multi-group radiation diffusion equations~\cite{hang2013convergence,morel1985synthetic}.
In each source iteration, the $G$ radiation diffusion equations are decoupled,
and each group radiation diffusion equation can be solved independently.
By considering the time advancing,
Algorithm~\ref{alg:mg-radiation-diff-sol} can be used to
describe the main solution procedure for multi-group radiation diffusion equations.

\begin{algorithm}[htbp]
\caption{Solution procedure for multi-group radiation diffusion equations}
\label{alg:mg-radiation-diff-sol}
\begin{algorithmic}[1]
\For {(Time stepping loop)}
   \For {(Nonlinear iteration loop)}
      \For {(Source iteration loop)}
         \State{Solve $G$ linearized multi-group radiation diffusion equations.}
      \EndFor
   \EndFor
\EndFor
\end{algorithmic}
\end{algorithm}

A Lagrangian method is used in the Lared-integration code and
the mesh moves with the motion of the fluid.
Consequently, the multi-group radiation diffusion equations
are discretized on two-dimensional deforming meshes.
A nine-point scheme is used to discretize the diffusion equations.

\subsubsection{Linear systems for test}

In the test, the number of {cells} is 7040.
The number of groups is $G=64$, and therefore,
{64 equations must be solved} in each source iteration.

\begin{figure}[htbp]
  \centering
  \includegraphics[width=0.7\linewidth]{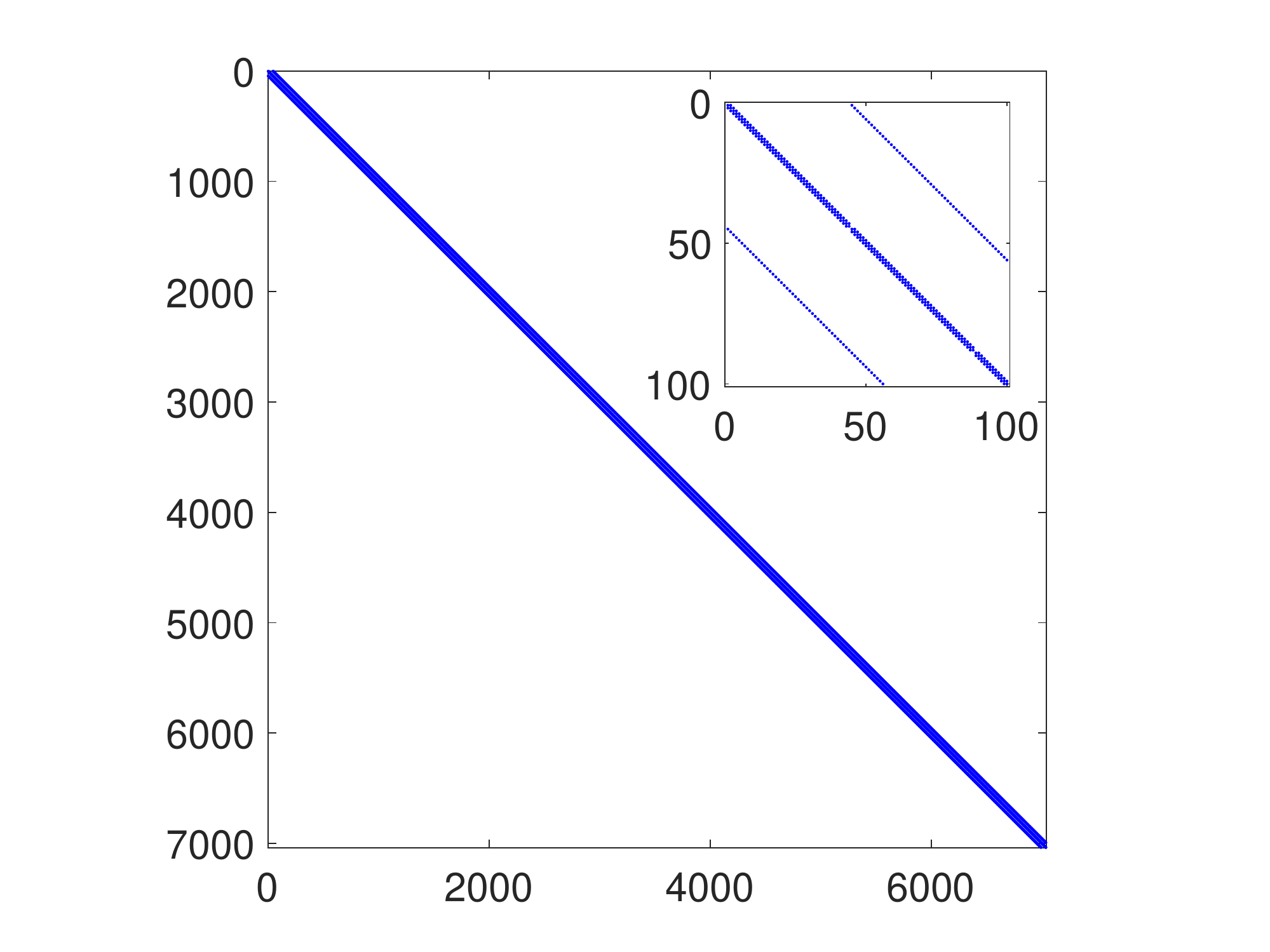}\\
  \caption{The sparsity pattern of the matrix of the multi-group radiation diffusion linear system.}
  \label{fig:spy1}
\end{figure}

Linear systems are tested for all groups, in a given source iteration of a given
nonlinear iteration, at time steps 1,000 and 5,000.
For notation, each linear system is denoted as
T$n_1$N$n_2$S$n_3$G$n_4$, which means the linear system
of group $n_4$ in the $n_3$th source iteration of $n_2$th nonlinear iteration at $n_1$th time step.
For example, T1000N0S0G0 represents the linear system of $0$-th group in $0$-th source iteration
of $0$-th nonlinear iteration at time step 1,000.
All the coefficient matrices of the linear system have the same sparsity pattern, shown in Fig.~\ref{fig:spy1}.

\subsubsection{{Test results}}

\begin{figure}[htbp]
  \centering
  \begin{minipage}[t]{0.8\linewidth}
  \includegraphics[width=1\linewidth]{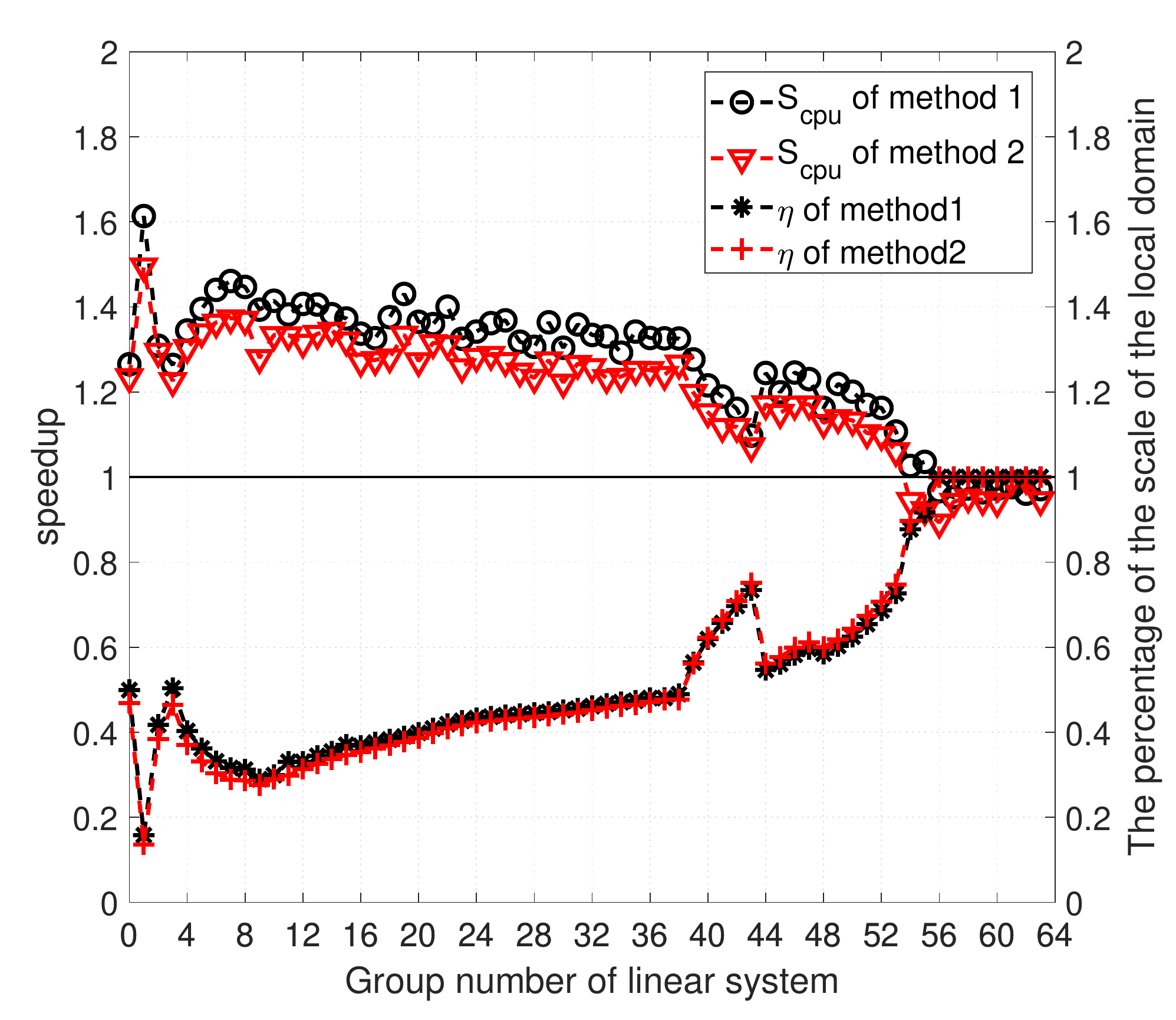}
  \end{minipage}
  \begin{minipage}[t]{0.8\linewidth}
  \includegraphics[width=1\linewidth]{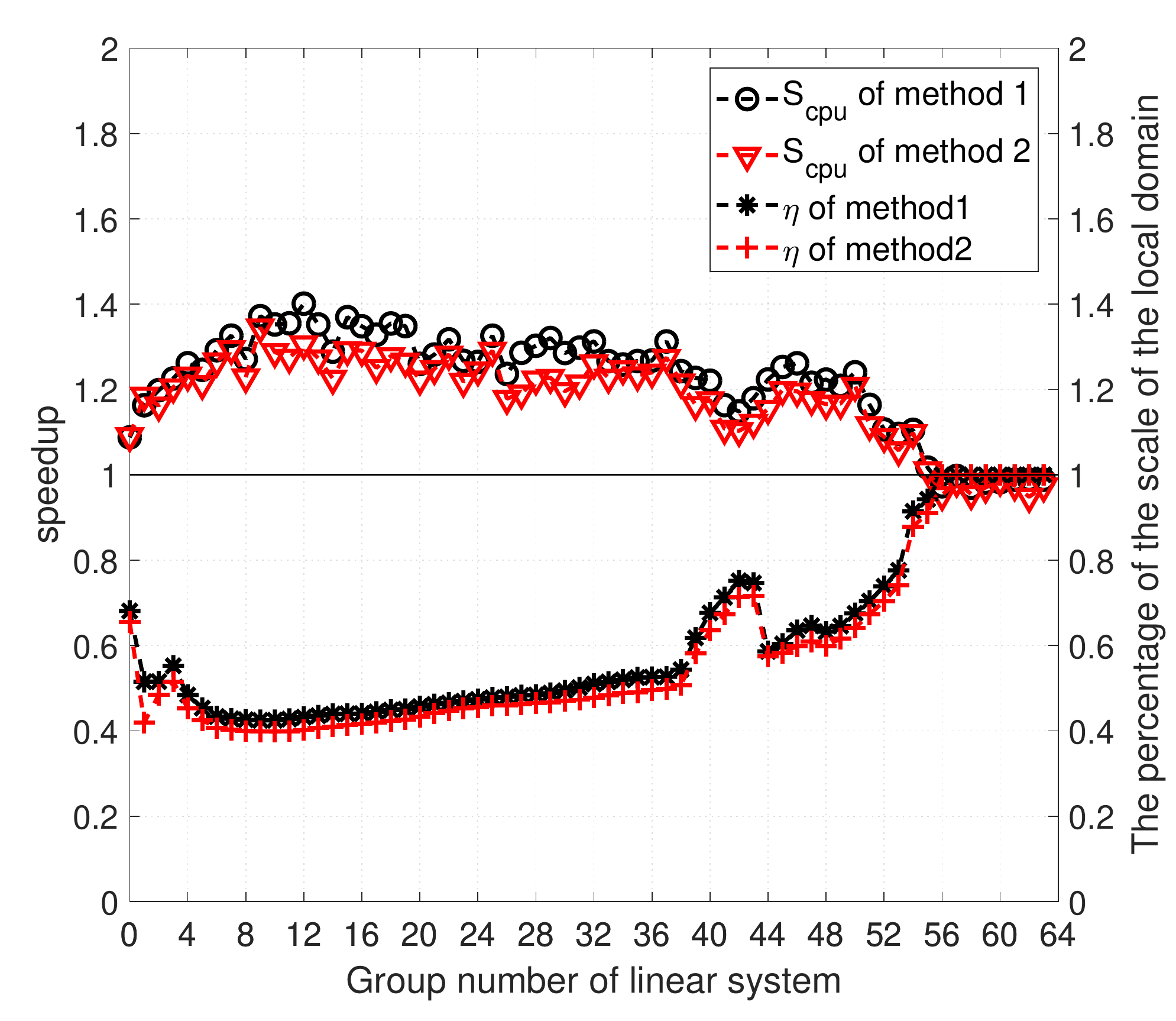}
  \end{minipage}
  \caption{The percentage of the scale of the subsystem and
  the speedup of the local {character-based} methods
  for solving multi-group radiation diffusion equations.
  The upper is the result of T1000N1S0G0-G63 at time step 1,000,
  the lower is the result of T5000N0S0G0-G63 at time step 5,000.}
  \label{fig:comparison1}
\end{figure}

We test a total of 128 linear systems on two time steps,
i.e., T1000N0S0G0-G63 and T5000N0S0G0-G63.
In this test, the parameter $\alpha$ in $\method1$ is set to $10^{-10}$,
and $E_{\max}$ in $\method2$ is set to 1. The convergence tolerance $\epsilon = 10^{-8}$.
The number of the Gauss-Seidel iteration is one in Algorithm~\ref{alg:local-character}.
The speedups and the percentages of {the} scale of the local domain are shown in Fig.~\ref{fig:comparison1}.

In Fig.~\ref{fig:comparison1}, at 1,000th time step, the maximal speedup of $\method1$ and $\method2$ for all the 64 equations are 1.61 and 1.50, respectively,
and the corresponding percentage of the scale of the local domain {subsystems} are 15.8\% and 13.6\%, respectively.
At 5,000th time step, the maximal speedup of $\method1$ and $\method2$ for all the 64 equations are 1.40 and 1.34, respectively, and the corresponding percentage of the scale of the local {domains} are 43.2\% and 39.8\%, respectively.

The comparison of the effectiveness of the methods
for some specific linear systems at 1,000th time step {is} given in \Cref{table:speedup-MG}.
This table shows that with the increase of the group number,
the speedup of the local {character-based} method decreases
because $\eta$ increases.
In particular, when $\eta$ approximates 1,
the time {costs} of both the local {character-based} methods are larger than that of $\method0$.
Further analysis {of} the choice of the parameters will be presented in \cref{sec:IOF}.

\begin{table}
\setlength{\tabcolsep}{3.0mm}
\caption{
\label{table:speedup-MG}
Speedup of $\method1$ and $\method2$ to $\method0$ for some groups of T1000N1S0.}
\begin{center}
\begin{tabular}{ccccccc}
\hline
G\#     & Method     & $\eta$ & $S_{\CPU}$ & $\CPU_{\total}$ & $\CPU_{\localconstruct}$ & $\CPU_{\localsolve}$  \\
\hline
        & $\method0$ & $-$    & 1.0        &  1.95E-2        & $-$                      & $-$                   \\
G1      & $\method1$ & 15.8\% & 1.61       &  1.21E-2        & 4.73E-3                  & 6.65E-3               \\
        & $\method2$ & 13.6\% & 1.50       &  1.30E-2        & 5.62E-3                  & 6.74E-3               \\
\hline
        & $\method0$ & $-$    & 1.0        &  2.16E-2        & $-$                      & $-$                   \\
G10     & $\method1$ & 29.9\% & 1.41       &  1.53E-2        & 5.01E-3                  & 9.64E-3               \\
        & $\method2$ & 29.0\% & 1.33       &  1.62E-2        & 5.68E-3	                & 9.89E-3               \\
\hline
        & $\method0$ & $-$    & 1.0        &  2.65E-2        & $-$                      & $-$                   \\
G41     & $\method1$ & 65.7\% & 1.19       &  2.23E-2        & 5.61E-3                  & 1.61E-2               \\
        & $\method2$ & 66.4\% & 1.12       &  2.37E-2        & 6.40E-3	                & 1.67E-2               \\
\hline
        & $\method0$ & $-$    & 1.0        &  3.25E-2        & $-$                      & $-$                   \\
G55     & $\method1$ & 91.7\% & 1.04       &  3.14E-2        & 6.18E-3                  & 2.45E-2               \\
        & $\method2$ & 93.5\% & 0.92       &  3.52E-2        & 6.74E-3                  & 2.78E-2               \\
\hline
        & $\method0$ & $-$    & 1.0        &  3.68E-2        & $-$                      & $-$                   \\
G63     & $\method1$ & 100\%  & 0.97       &  3.79E-2        & 6.24E-3                  & 3.11E-2               \\
        & $\method2$ & 100\%  & 0.95       &  3.89E-2        & 6.79E-3                  & 3.15E-2               \\
\hline
\end{tabular}
\end{center}
\end{table}

In multi-group radiation energy model,
the high energy groups represent the energy of high frequency photons.
Usually, the number of high frequency photons is small.
This makes the energy of high frequency photons low.
Particularly, for multi-group radiation diffusion equations,
the solution of high energy groups is small in the computing domain,
and the local character is not strong. This makes the method
ineffective for high energy groups.
At the same time, with time advancing, the temperature (representing energy) increases
in the simulation domain for multi-group radiation diffusion equations
(also for three temperature energy equations),
the local character becomes weak. Therefore, the performance of
the local character-based method deteriorates at large time steps.

\subsection{Test of three temperature energy equations}
\label{sub:result2}

When only one group is used in the multi-group radiation diffusion equations,
the multi-group radiation diffusion equations coupled to electron-ion equations
are reduced to the three temperature (3T) energy equations,
which can be given as follows:
\begin{eqnarray}
\label{eq:3T}
\left\{
\begin{array}{lll}
C_{ve}{\frac{\partial T_e }{\partial t}} - \frac{1}{\rho} \nabla \cdot (K_e \nabla T_e) &=& \omega_{ei} (T_i - T_e) +
\omega _{er}(T_r - T_e) \\
C_{vi}{\frac{\partial T_i}{\partial t}} - \frac{1}{\rho}\nabla \cdot (K_i \nabla T_i) &=& \omega_{ei}(T_e - T_i) \\
C_{vr}{\frac{\partial T_r}{\partial t}} - \frac{1}{\rho}\nabla \cdot (K_r \nabla T_r) &=& \omega_{er}(T_e - T_r)
\end{array}
\right.
\end{eqnarray}
In the equations,
\begin{itemize}
\item $T_e$, $T_i$, and $T_r$ are electron temperature, ion temperature, and photon temperature, respectively.
\item $C_{ve}$,  $C_{vi}$, and $C_{vr}$ are electron heat capacity, ion heat capacity,
      and photon heat capacity, respectively.
\item $\rho$ is the material density.
\item $\omega_{ei}$ and $\omega_{er}$ are the energy exchanging coefficients between {electrons and ions,
      and that between electrons and photons}, respectively.
\item $K_{e}$, $K_{i}$, and $K_{r}$ are electron heat conduction coefficients,
      ion heat conduction coefficients, and photon heat conduction coefficients, respectively.
\end{itemize}

Similar discretization methods as for multi-group radiation diffusion equations are used.
More details of the 3T equations are presented in \cite{an2009choosing,zeyao2004parallel}.

\subsubsection{Linear systems for test}

\begin{figure}[htbp]
  \centering
  \includegraphics[width=0.75\linewidth]{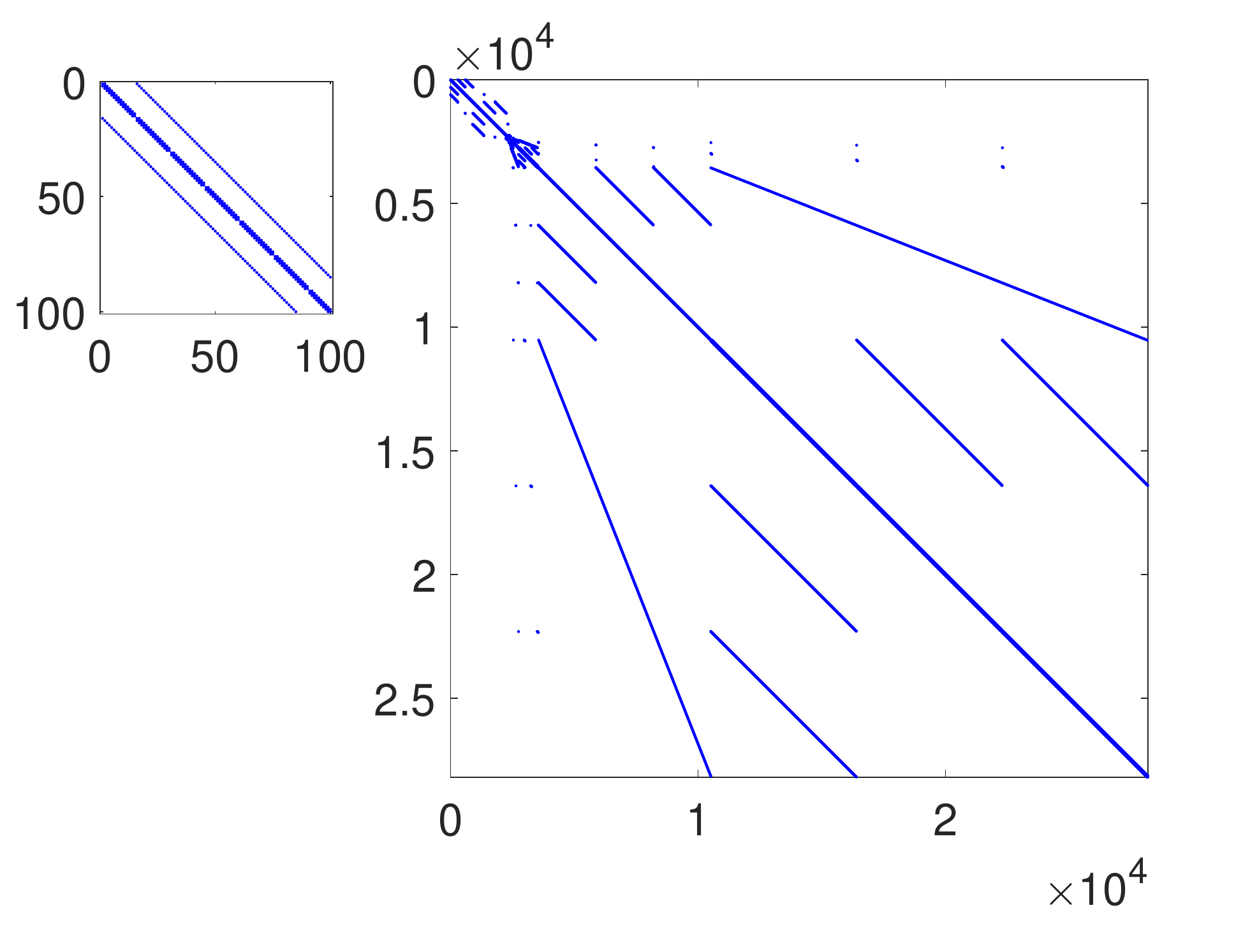}\\
  \caption{The sparsity pattern of the matrix in linear system of 3T equations.}\label{fig:3T}
\end{figure}

The test set in this case consists of eight linear systems, which are extracted from the simulation
{every 10,000 time steps, from the 20,000th time step to the 90,000th time step.}
After being discretized, the variables of the 3T energy equations are organized as follows in the linear system
\begin{eqnarray*}
T = (T_1^e, T_1^i, T_1^r, \ldots, T_M^e, T_M^i, T_M^r)^{\small T},
\end{eqnarray*}
where $M$ is the number of the cells. In this case, $M = 9,396$,
and the total number of the unknowns in each linear system is 28,188.
For convenience, we denote the linear systems in the test set as $t_{n}$,
which represents the linear system at the $n$-th time step.
All the coefficient matrices of the linear systems
have the same sparsity pattern as is shown in Fig.~\ref{fig:3T}.

\subsubsection{Test result}

\begin{figure}[htb]
  \centering
  \includegraphics[width=0.8\linewidth]{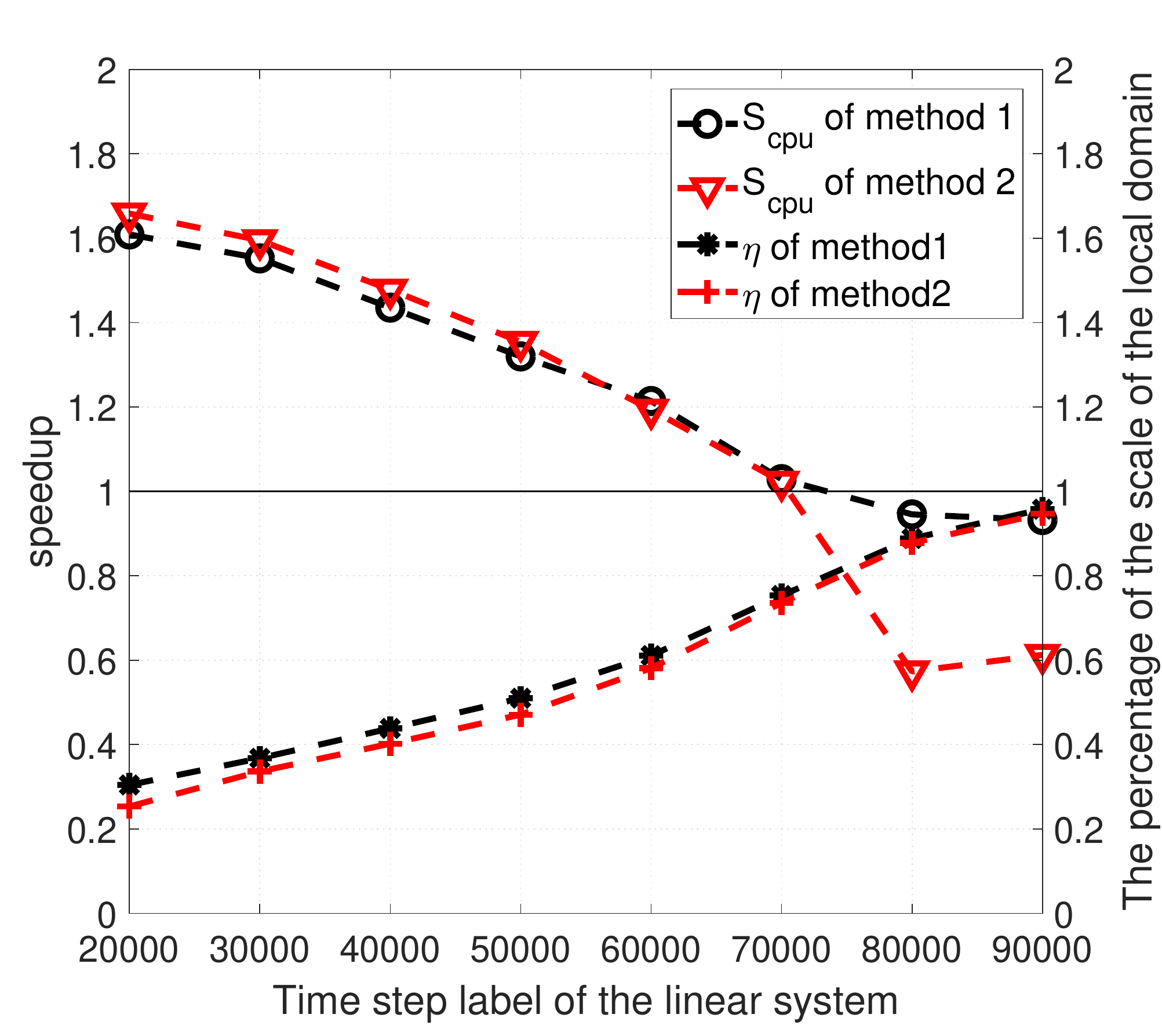}\\
  \caption{
  The percentage of the scale of the subsystem and
  the speedup of the local {character-based} methods
  for solving 3T energy equations.}
  \label{fig:3Tspeedup}
\end{figure}

In this test, all the parameters are the same as for the test of the multi-group radiation diffusion equations.
The speedups and the percentages of the scale of the local domain are shown in Fig.~\ref{fig:3Tspeedup}.
From this figure, one can see that, for most of the tested linear systems,
$\method1$ and $\method2$ are more efficient than $\method0$.
The maximal speedup of $\method1$ to $\method0$ is 1.61,
and the corresponding percentage of the scale of the local domain is about 30.4\%.
The maximal speedup of $\method2$ to $\method0$ is about 1.65,
and the corresponding percentage of the scale of the local domain is about 25.3\%.

The detailed results for four linear systems of 3T energy equations are given in \Cref{table:speedup-3T}.
From this table, one can see that for the linear equations at time step 20,000, 40,000 and 60,000,
the speedups of $\method1$ and $\method2$ are larger than 1.0.
For the linear system at time step 80,000, the speedup is
less than one for both $\method1$ and $\method2$.
That is, the solution time of $\method1$ and $\method2$ is longer than $\method0$.
The parameters in the methods have some influence on the performance,
and this will be studied further in {the} next subsection.

\begin{table}
\setlength{\tabcolsep}{3.0mm}
\caption{
\label{table:speedup-3T}
Speedup of $\method1$ and $\method2$ to $\method0$ for 3T energy equations on four time steps.}
\begin{center}
\begin{tabular}{ccccccc}
\hline
$t_{\#}$   & Method     & $\eta$ & $S_{\CPU}$ & $\CPU_{\total}$ & $\CPU_{\localconstruct}$ & $\CPU_{\localsolve}$  \\
\hline
            & $\method0$ & $-$    & 1.0        &  4.25E-2        & $-$                      & $-$                   \\
$t_{20000}$ & $\method1$ & 30.4\% & 1.61       &  2.64E-2        & 9.61E-3                  & 1.45E-2               \\
            & $\method2$ & 25.3\% & 1.66       &  2.56E-2        & 9.88E-3                  & 1.35E-2               \\
\hline
            & $\method0$ & $-$    & 1.0        &  4.85E-2        & $-$                      & $-$                   \\
$t_{40000}$ & $\method1$ & 43.8\% & 1.44       &  3.38E-2        & 1.10E-2                  & 2.05E-2               \\
            & $\method2$ & 40.2\% & 1.48       &  3.28E-2        & 1.11E-2	                & 1.94E-2               \\
\hline
            & $\method0$ & $-$    & 1.0        &  4.85E-2        & $-$                      & $-$                   \\
$t_{60000}$ & $\method1$ & 61.1\% & 1.21       &  3.99E-2        & 1.17E-2                  & 2.59E-2               \\
            & $\method2$ & 58.1\% & 1.19       &  4.06E-2        & 1.25E-2	                & 2.57E-2               \\
\hline
            & $\method0$ & $-$    & 1.0        &  5.77E-2        & $-$                      & $-$                   \\
$t_{80000}$ & $\method1$ & 88.8\% & 0.95       &  6.10E-2        & 1.50E-2                  & 4.34E-2               \\
            & $\method2$ & 87.8\% & 0.68       &  1.01E-1        & 1.56E-2                  & 4.10E-2               \\
\hline
\end{tabular}
\end{center}
\end{table}

\subsection{Influences of the parameters}
\label{sec:IOF}

In both Algorithm~\ref{alg:gradient-local-domain} and Algorithm~\ref{alg:residual-local-domain},
the parameters should be set for specific applications.
In the previous tests, the parameters are fixed.
In this subsection, the influence of the parameters
for the algorithms will be analyzed.

\subsubsection{The parameter $\alpha$ in $\method1$}

In $\method1$, $\alpha$ is used to construct the local domain $\Omega_{\local}$.
Therefore, $\alpha$ has a significant influence on the effectiveness of the method.
When $\alpha$ is too small, the scale of the subsystem (measured by $\eta$,
the percentage of the scale of
the local domain subsystem to the scale of the global system) will be large.
This will deteriorate the performance of $\method1$.
Conversely, when $\alpha$ is too large, $\eta$ will be small,
and the local domain $\Omega_{\local}$ may not be large enough to
capture the physical information being propagated forward in time.
This will further decrease the performance of the local {character-based} method.
In particular, if $\alpha=0$, then $\Omega_{\local} = \Omega$, the whole domain;
if $\alpha=1$, then $\Omega_{\local} = \emptyset$.
For a linear system, one way to determine the optimal value of $\alpha$ is
\begin{eqnarray*}
\alpha_{\opt} = \arg\min_{\alpha} \{\alpha : \| b - A \tilde{x}_{\alpha} \|_2 \leq \epsilon \|b\|_2 \},
\end{eqnarray*}
where $\tilde{x}_{\alpha}$ is the solution obtained by the local {character-based} method,
and the local domain is constructed with the the parameter $\alpha$.

To show the influence of the parameter $\alpha$ on the method,
two linear systems of multi-group radiation diffusion equations,
T1000N1S0G0 and T1000N1S0G55, are selected.
The value of $\alpha$ is set to 12 different orders of magnitude values, varying from 1 to $10^{-11}$.

\begin{figure}[htbp]
\centering
  \begin{minipage}[t]{0.65\linewidth}
  \includegraphics[width=1\linewidth]{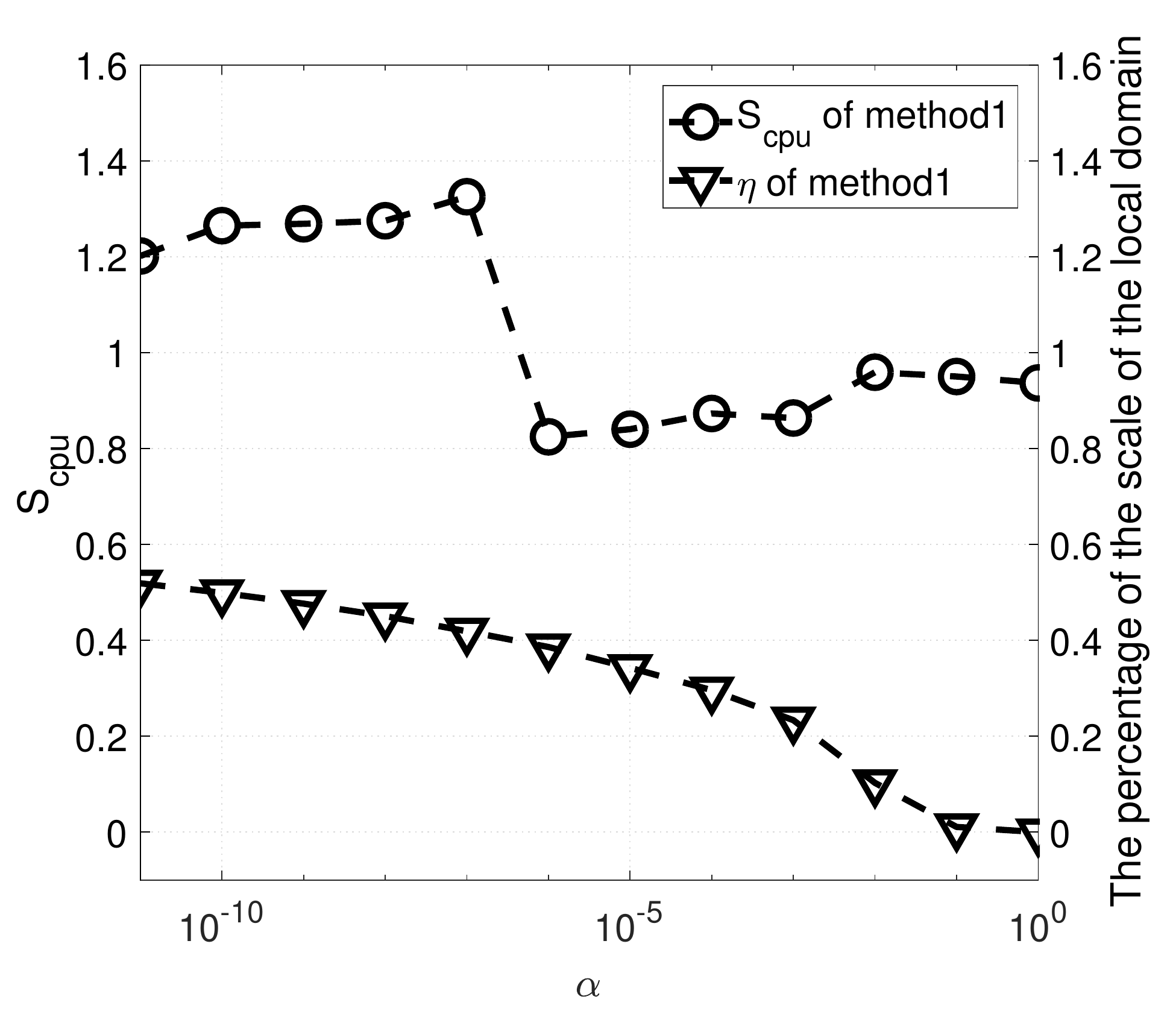}
  \end{minipage}
 \begin{minipage}[t]{0.65\linewidth}
  \includegraphics[width=1\linewidth]{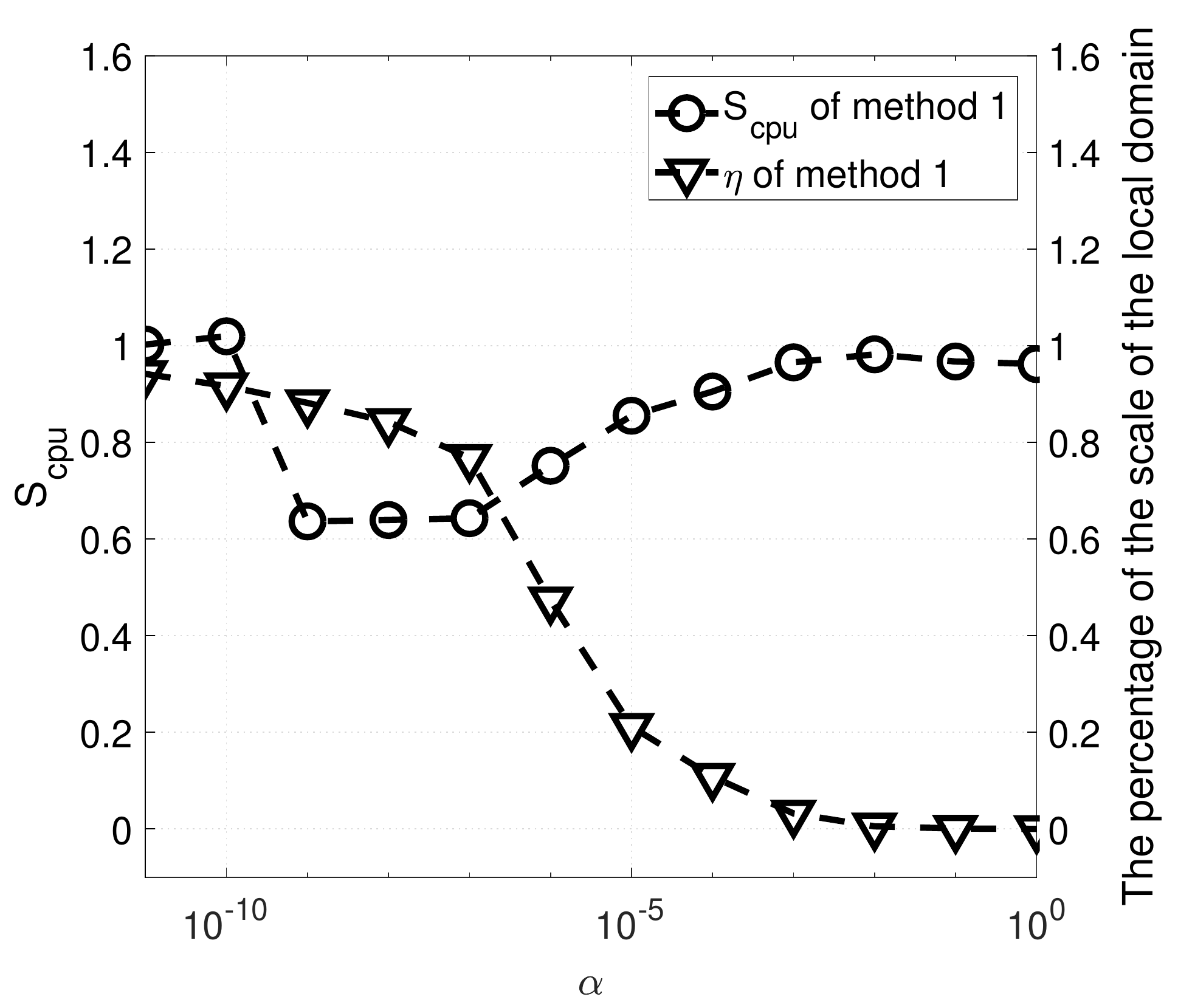}
  \end{minipage}
 \caption{The percentage of scale of the local domain and speedup of $\method1$ with different $\alpha$.
          The upper sub-figure is the result of T1000N1S0G0,
          the lower sub-figure is the result of T1000N1S0G55.}
 \label{fig:alpha}
\end{figure}

Fig.~\ref{fig:alpha} shows the speedups of the local character-based method and the percentage of the scale of the local domain for both cases.
The results show that $\eta = \eta(\alpha)$ is monotonic decreasing,
which is consistent with the construction of the local character-based method.
For the tested two linear systems, the $\alpha$ values corresponding to the maximal speedup are
$10^{-7}$ and $10^{-10}$, respectively, and the maximal speedups are 1.38 and 1.1, respectively.
Particularly, for the linear system T1000N1S0G55, this figure shows that the speedups are lower than 1
for most of the cases. This indicates that $\method1$ is not a good candidate for solving
this linear equation because the local character of this linear system is not so strong.

\subsubsection{The parameter $E_{\max}$ in $\method2$}

\begin{figure}[htbp]
\centering
  \begin{minipage}[t]{0.65\linewidth}
  \includegraphics[width=1\linewidth]{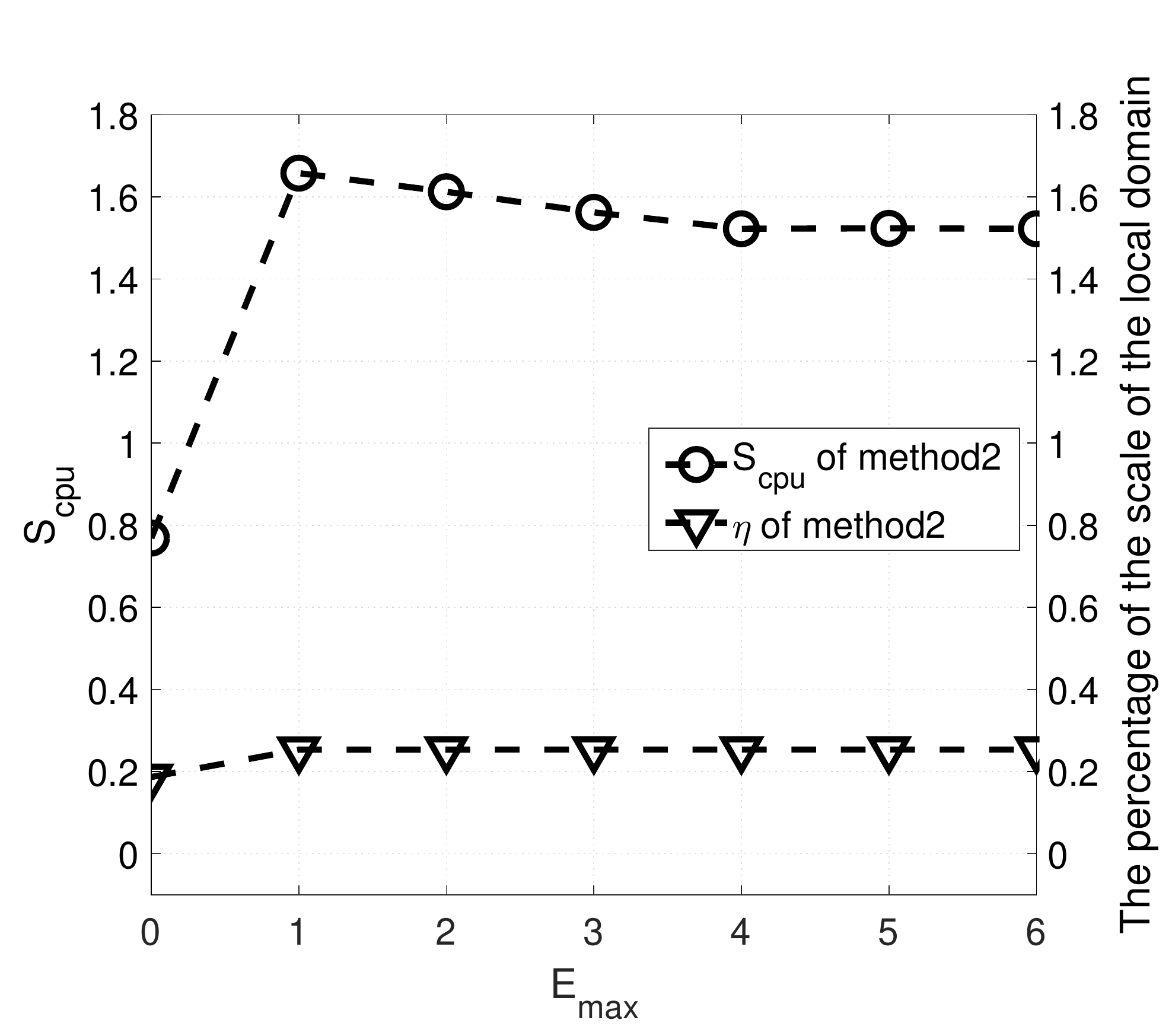}
  \end{minipage}
 \begin{minipage}[t]{0.65\linewidth}
  \includegraphics[width=1\linewidth]{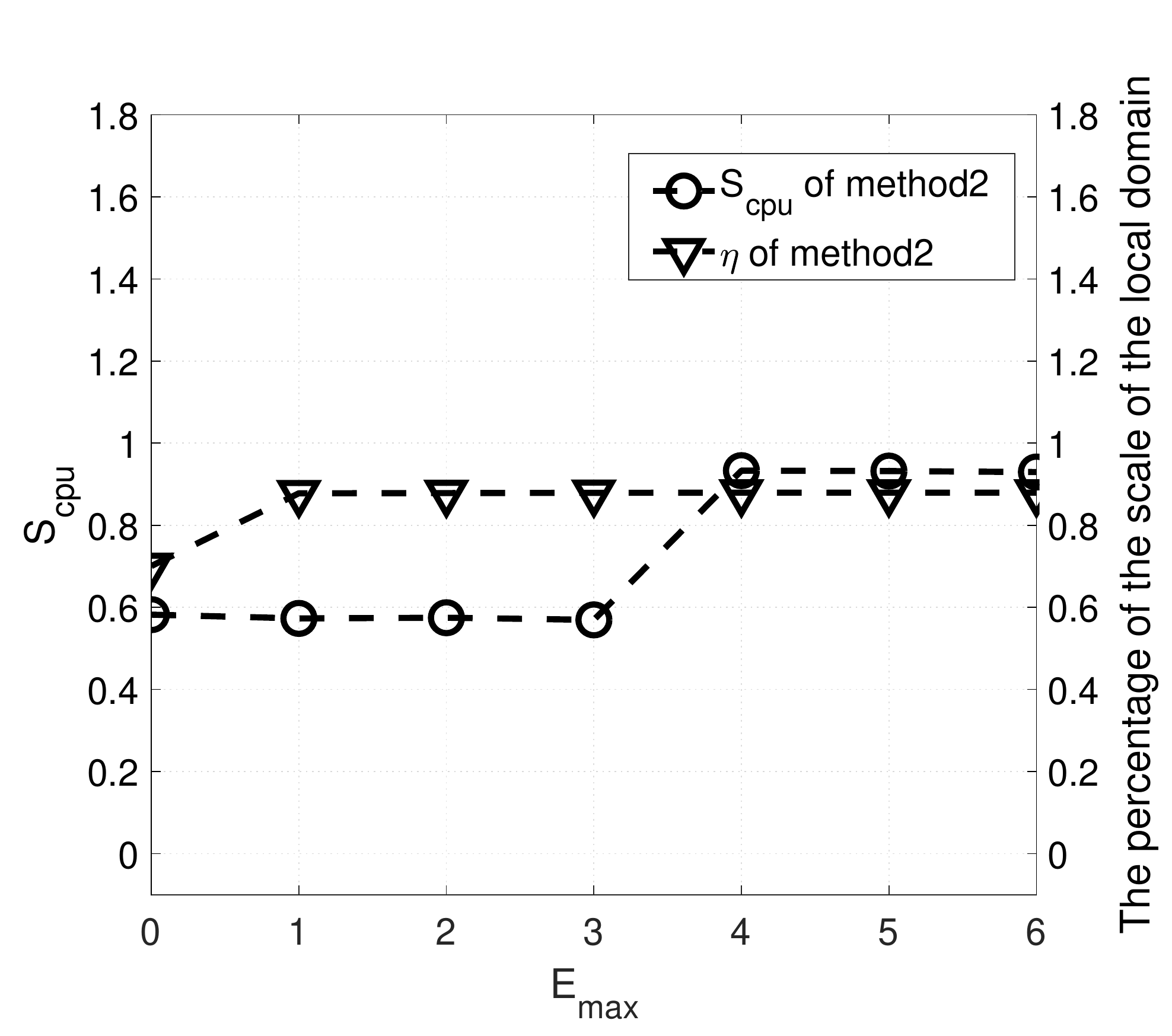}
  \end{minipage}
 \caption{The percentage of scale of the local domain and speedup of $\method2$ with different $E_{\max}$.
          The upper sub-figure is the result of $t_{20000}$,
          the lower sub-figure is the result of $t_{80000}$.}
 \label{fig:3ts}
\end{figure}

In $\method2$, $E_{\max}$ is used in the second stage of the local domain construction method,
and it represents the maximal expanding iterations for the local domain $\Omega_{\local}$.
Theoretically, a large enough value for $E_{\max}$ can be set
so that $\Omega_{\local}$ {cannot} be expanded {anymore}.
However, the larger that $E_{\max}$ is,
the more expensive constructing the local character-based method will be.
An improper $E_{\max}$ may result in a loss of speedup of the method.
To show the specific influence of the parameter $E_{\max}$,
two linear systems of 3T energy equations, $t_{20000}$ and $t_{80000}$,
are tested further with different $E_{\max}$.
In the test, the value of $E_{\max}$ varies from 0 to 6.
The test results are shown in Fig.~\ref{fig:3ts}.

From Fig.~\ref{fig:3ts}, one can see that, for linear system $t_{20000}$, $\eta$ does not have clear increase
when $E_{\max}$ increases from 1 to 6, and the maximal speedup is 1.66 when $E_{\max} = 1$.
This means that $E_{\max}=1$ is the optimal value for solving this linear system.
For linear system $t_{80000}$, again the value of $\eta$ does not increase clearly
when $E_{\max}$ increases from 1 to 6. However, the speedup has a dramatic increase
when $E_{\max}$ increases from 3 to 4, and the optimal $E_{\max}$ value is 4
for solving this linear system. However, for this linear system,
the speedups are less than 1 for all cases.
This means that the local {character-based} method is not good at
solving the linear system $t_{80000}$ because its local character is not so strong.

\section{Summary and remarks}
\label{sec:summary}

The solution of radiation diffusion equations is very important
in the simulation of inertial confinement fusion (ICF) and astrophysics.
The solution of radiation diffusion equations concerns many linear equation sequences
with the total number of linear equations is about $\mathcal{O}(10^4)$.
The radiation diffusion equations are {time-dependent} nonlinear systems
and the solution at the previous time step or nonlinear iterate is usually
used as the initial iterate for solving the linear equations.
Furthermore, in most of the simulation {period},
the solution only varies greatly
in some local domain.

In this paper, a local {character-based} method is proposed
to solve the linear systems in this {application}.
In the local {character-based} method, {first,} a local domain is constructed.
It is expected that on the local domain the solution varies greatly
from the initial iterate to the last iterate solution.
Then the subsystem on the local domain is solved.
At last, the global linear system is solved
with the initial iterate obtained by solving the local domain system.
Two methods are proposed to construct the local domain.
One is based on the spatial gradient{,} and the other is based on the initial residual.
Therefore, two kinds of specific methods are constructed.

Numerical results show that the proposed local character-based methods
can successfully grasp the local characteristic domain.
A two-dimensional heat conduction equation is tested and analyzed.
Also, the linear systems from real applications of multi-group radiation diffusion equations and 3T energy equations are tested.
Compared to solving the linear system directly by some method, such as preconditioned Krylov methods,
the solution time can be reduced by 38\% and 40\% respectively by using the two local {character-based} methods.

The performance of the local {character-based} method
depends strongly on the percentage of the scale of the local domain subsystem to
the scale of the global system.
Numerical results for multi-group radiation diffusion equations and
3T energy equations show that if the percentage is relatively small,
then the local {character-based} method may perform better than the typical method,
such as AMG preconditioned GMRES method.
Otherwise, the local {character-based} method may perform worse than the typical method.
Therefore, in the simulation of real applications,
if the temperature is not high in the whole simulation domain,
then the local character is strong{,} and the method proposed in this paper
will be efficient. In the later simulation {period},
with the increase of the temperature in the whole domain,
the local character may not be so strong, the typical method should be used.
In {the} future, we will consider {constructing the} preconditioning method based on
the local {character-based} method.

%
%


%
%


\end{document}